\documentclass[11pt]{amsart}

\renewcommand{\bar}{\overline}
\newcommand{\eps}{\varepsilon}
\newcommand{\pa}{\partial}
\renewcommand{\phi}{\varphi}

\newcounter{hours}\newcounter{minutes}

\newcommand{\ka}{K\"ahler }

\newcommand{\C}{{\mathbb C}}

\newcommand{\R}{{\mathbb R}}

\newcommand{\frk}[1]{{\mathfrak{#1}}}

\newcommand{\bb}{{\frac{\sqrt{-1}}{2\pi}}}

\renewcommand{\epsilon}{\varepsilon}

\usepackage{amsmath,amsthm,amscd} 
\usepackage{calc}         

\title 
[]{On the Weil-Petersson volume and the first Chern Class of the moduli
space of Calabi-Yau manifolds}
\author{Zhiqin Lu and Xiaofeng Sun}
\date{\today}

\address[Zhiqin Lu and Xiaofeng Sun] 
{Department of Mathematics\\
University of California, Irvine\\ 
Irvine, CA 92697}
\email[Zhiqin Lu]{zlu@math.uci.edu}
\email[Xiaofeng Sun]{xsun@math.uci.edu}
\thanks{The first author is supported by NSF Career Award
DMS 0347033 and an Alfred P. Sloan Fellowship. The second author is
supported by NSF grant DMS 0202508.}

\newtheorem{theorem}{Theorem}[section] 

\newtheorem{lemma}{Lemma}[section]
\newtheorem{cor}{Corollary}[section]
\newtheorem{prop}{Proposition}[section]
\newtheorem{claim}{Claim}
\newtheorem{definition}{Definition}[section]

\theoremstyle{remark}
\newtheorem{rem}{Remark}[section]  
\newtheorem*{qu}{\bf Conjecture}   

\begin{document}     
\maketitle  



\numberwithin{equation}{section}          

\tableofcontents

\section{Introduction}
In this paper, we continue our study of the Weil-Petersson 
geometry as in the previous paper~\cite{LS-1}, in which we
have proved  the boundedness of the Weil-Petersson volume,
among the other results. The main results of this paper are that
the volume and the integrations of Ricci curvature
of the Weil-Petersson metric on 
the moduli space are rational numbers. In 
particular,
the Ricci curvature defines the first Chern class of the 
moduli space in the
sense of Mumford~\cite{Mumford1}.

It was a classical result of Mumford~\cite{Mumford1} that
for a noncompact \ka manifold $M$ with $\bar M$ being
a smooth compactification of $M$ and $\bar M\backslash M$
being  a divisor
$D$ of normal crossings, and for  any Hermitian bundle
$(E,h)$ over
$M$, one can define the Chern classes $c_k(E)$  provided the metric $h$
is ``good'' defined by Mumford ~\cite[Section 1]{Mumford1}.
Roughly speaking, a metric is ``good'' if 
the metric matrix have log bound, and the local connection form and the curvature have 
 Poincar\'e  type growth.
It was verified that  
the natural bundles over locally Hermitian symmetric spaces
are ``good'' (cf. ~\cite{Mumford1}). 
For the moduli space of curves of genus greater than or equal to
$2$,  the metric induced by the Weil-Petersson metric on the determinant bundle of the log extension of the cotangent bundle is good~\cite{Tra1}.
However, for the 
moduli space of polarized Calabi-Yau manifolds, 
it is not clear that the
Weil-Petersson metric or the volume form of the
Weil-Petersson metric is ``good''. By ~\cite{BCOV1},
the Weil-Petersson potential is related to the
analytic torsion of the moduli space. While the
Hessian of the torsion is known to be related to
the Weil-Petersson metric and the generalized Hodge
metric(\cite{FL1}), it is not easy to find the
asymptotic behavior of the BOCV torsion itself.
Thus we can not use the theorem of Mumford directly
to prove that the integrations are rational numbers.

In this paper, we avoided using the BCOV torsion by
the careful analysis of the asymptotic behavior of the
Hodge bundles  at infinity.
Using the Nilpotent Orbit theorem of Schmid~\cite{Schmid},
we can
give another explicit (local) representation of the Weil-Petersson
potential. The potential has the following properties:
first, as a potential of a \ka metric, 
it  must be plur-subharmonic. Next, by the Strominger
formula, the Ricci curvature of the
Weil-Petersson metric  is lowerly bounded. Thus
  the volume
form of the Weil-Petersson metric is also plur-subharmonic. 
Our analysis
relies heavily on the above two properties of the potential
of the Weil-Petersson metric.

For the noncompact manifold $M$ defined above, we
can  define a \ka metric, called the 
Poincar\'e metric on $M$ such that on each 
Euclidean neighborhood of $D$, the metric is
asymptotically Poincar\'e (See \S ~\ref{cut} for the 
precise definition). 
The metric is not canonically chosen so it doesn't reflect
the geometry of the manifold $M$. However,
it is complete and the volume is finite, and
its curvatures are bounded. In practice, we use the 
Poincar\'e metric to bound the other intrinsically
defined metrics.

Let $M$ be a Weil-Petersson variety (see \S 2 for 
the definition).
It is not clear  whether we can construct a \ka metric
with nonpositive sectional curvature on the Weil-Petersson 
variety. However, if we only require the nonpositivity of the
bisectional curvature, then 
 in
~\cite{Lu3} and ~\cite{Lu5}, the first author defined such a
metric, which we called the Hodge metric. In addition
to the nonpositivity of the bisectional curvature of the Hodge
metric, the holomorphic sectional curvature  and the Ricci 
curvature of Hodge metric are negative 
and bounded away from
zero. By Schwarz-Yau~\cite{Y3} lemma, the Hodge metric is
bounded by the Poincar\'e metric defined above.

Using the comparison of the Hodge and the Poincar\'e metrics,
in the previous paper \cite{LS-1}, we have proved that the
Weil-Petersson volume and the Hodge volume are all finite. 
By the definition of the Weil-Petersson metric and the
Hodge metric, if a Weil-Petersson variety {\sl were}
compact, then the  volume with respect to the
Weil-Petersson metric and the Hodge metric could have been
rational numbers because of the Gauss-Bonnet
Theorem. We thus
conjectured that the volume of both metrics are in fact
rational numbers, even though the Weil-Petersson varieties
are  more likely to be  noncompact.

In this paper, we verified the conjecture by controlling
the growth of the potential of the Weil-Petersson metric at infinity.
In order to get
the estimates we need, we have to define a special kind of 
cut-off functions. In general, if a cut-off function is $1$ at
the origin and if it is supported within a ball of radius
$r$, then its second derivatives are of the order $1/r^2$.
We can do a little bit better for the Hessian of the cut-off 
function on $\R^2$, because $\R^2$ is an example of 
parabolic manifold defined by P.~Li~\cite{Li}.
In fact, the order of the Hessian of the cut-off function
is of the order $\frac{1}{r^2(\log 1/r)^2}$.
This observation is important in our proof. By using the cut-off 
function and the convexity of the Weil-Petersson potential,
we can prove

\begin{theorem}\label{thm11}
Let $(M,\omega_{WP})$ be a Weil-Petersson variety of
dimension $m$. Then its volume
\[
\int_{M}\omega^m_{WP}
\]
is a rational number.
\end{theorem}

The volume form of the Weil-Petersson metric has its own
convexity by the formula of Strominger. 
However, in this case, the (local) volume forms are
not integrable with respect to the Poincar\'e
metric. Special care must be taken in order
to get the similar result as in the volume case. In order
to do that, we defined the degeneration order of the
volume form along each hypersurface of the divisor
$\bar{M}\setminus {M}$ and then twisted the
extension of the anti-canonical bundle of  $M$
(See Definition~\ref{def666} for details). Using this, we
can prove

\begin{theorem} \label{thm12}
Let $X \subset  M$ be a subvariety of dimension
$q$. Then for non-negative integers $k$ and $l$ with
$k+l=q$,     
\[
\int_{X}(Ric(\omega_{WP}))^k\wedge\omega_{WP}^l
\]
are rational numbers.
\end{theorem}

Obviously, Theorem~\ref{thm12} implies Theorem~\ref{thm11}.
On the other side, if $X,M$ in Theorem~\ref{thm12}
are smooth,  then we have

\begin{cor} Using the above notations, we have
\[
(c_1(\bar K_X^{-1})+\sum \mu_i Y_i)^k\cap
c_1(\bar F^n)^l=
\int_{X}(Ric(\omega_{WP}))^k\wedge\omega_{WP}^l,
\]
where 
$\bar F^n$ and 
$\bar K_X^{-1}$ are the Hodge extensions of  the 
Hodge bundles $F^n$ and the anti-canonical bundle
$K_X^{-1}$ of $X$;
$\sum Y_i=Y$ is the divisor $\bar X-X$;
$Y_i$ are irreducible components of $Y$; and
$\mu_i$ are the degeneration orders of the Weil-Petersson 
metric along $Y_i$.
\end{cor} 

\begin{rem}
In the above corollary, the righthand side is intrinsically
defined. Thus one will get some information of the 
divisors once the topology of the compact manifold
$\bar X$ is known. In particular, by using this we
can get the information of the monodromy group
assuming the moduli space is $CP^1$ minus three points.
Such a moduli space is of interest in Mirror Symmetry
(cf. Doran-Morgan~\cite{DM}).
\end{rem}

Two important papers in the direction of this paper
have drawn our attentions. One is the recent survey paper of
Todorov~\cite{To1}, which gives a complete summary of the recent progress
in  the subject. The other one is by Schumacher~\cite{Sc1},
in which the author computed the curvature of the
Weil-Petersson metric of K\"ahler-Einstein manifolds, using the
idea of Siu~\cite{Siu2} of horizontal liftings.

We are interested in the volume and the integrations of the curvature
because it defines some kind of invariants and using that, we wish to 
classify polarized Calabi-Yau manifolds and tell the monodromy
of the moduli space of Calabi-Yau manifolds~\footnote{In 
fact, our proof strongly hints the relations between
the monodromy operators and the rational numbers we defined.
In the case that the moduli space is of one-dimensional,
an explicit formula can be written down.}. These
questions are  very important in Mirror Symmetry and we shall
study them in a subsequent paper.

The organization of this paper is as follows: in \S 2, we give 
the definition of the Weil-Petersson geometry and some basic estimates;
in \S 3, we define the Poincar\'e metric and the cut-off function we need
for the rest of the paper; in \S 4, we write out the preferred extension
of the Hodge bundles defined by the Nilpotent Orbit 
Theorem
explicitly. The main part of the paper is \S 5 and \S 6, where
we prove Theorem~\ref{thm11} and Theorem~\ref{thm12}.

After finishing this paper, we were informed by A. Todorov that in~\cite{TD},
he proved the rationality of the volume of the moduli space independently.

{\bf Acknowledgment.}
The authors would like to thank  P. Li, D. Phong,
R. Schoen, and G. Tian for
their interest in the work. Particular thanks to K. Liu
for his many suggestions and encouragement during the preparation
of this paper.

\section{Preliminaries}\label{pre}

In this section, we  give the definition of the Weil-Petersson
geometry and Weil-Petersson variety, first appeared in~\cite{LS-1}.
Examples of Weil-Petersson varieties are the moduli spaces of
polarized Calabi-Yau manifolds.

\begin{definition} A Weil-Petersson variety is a \ka
orbifold
$M$ with the orbifold metric
$\omega_{WP}$ such that:
\begin{enumerate}
\item The universal covering space $\widetilde M$ is a smooth manifold. 
There is a natural immersion $\widetilde{M}\rightarrow D$ from 
$\widetilde M$ to the
classifying space $D$ (cf. ~\cite{Gr}) such that the image of 
$\widetilde M$ is a
horizontal slice of
$D$. The Hodge bundles $F^n\subset\cdots\subset  F^0$ are defined as the
pull-back of the  tautological bundles of $D$.
Furthermore, we have the natural identification
$TM=F^{n-1}/F^n$, where $TM$ is the holomorphic
tangent bundle of $M$;
\item $\omega_{WP}$ is the curvature of the bundle $F^n$. It
is positive definite and thus defines a \ka metric on
${M}$ and is called the Weil-Petersson metric;
\item $M$ is quasi-projective and $F^n$ extends to an ample line
bundle over  the compactification $\bar M$ of $M$;
\item After passing to a finite covering
and after desingularization, in a neighborhood of the infinity,
$M$ can be written as
\[
\Delta^{n-k}\times (\Delta^*)^{k},
\]
where $\Delta$ is the unit disk and $\Delta^*$ is the
punctured unit disk. Let $\Omega$ be a local section of
$F^p$ in the neighborhood, then locally, $\Omega$ can be
 written as
\[
\Omega=e^{\sqrt{-1}(N_1\log\frac{1}{z_1}+
\cdots+N_k\log\frac{1}{z_k}
)}A(z_1,\cdots, z_n),  
\]
where $N_1,\cdots,N_k$ are nilpotent operators
and $A$ is a vector-valued holomorphic function of
$z_1,\cdots,z_n$. Furthermore, all the local sections 
of all $F^p$, $p=1,\cdots,n$ 
satisfy the properties in the Nilpotent Orbit
theorem of 
Schmid~\cite{Schmid}.
\end{enumerate}
A Weil-Petersson subvariety $M_1$ is itself a Weil-Petersson variety. Moreover, it is  a subvariety of a Weil-Petersson
variety $M$ such that if $\bar M$ is a compactification of $M$, then 
$\bar{M_1}$, the closure of $M_1$ in $\bar M$ gives a compactification of
$M_1$ as a projective variety. 

The Weil-Petersson geometry is the geometry of the pair
$(M,\omega_{WP})$. 
\end{definition}
 
\begin{rem}
Moduli spaces of polarized Calabi-Yau manifolds are examples 
of Weil-Petersson varieties. In fact, for moduli spaces of 
polarized Calabi-Yau manifolds,
the first property of the above definition is  the transversality
property
of Griffiths'~\cite{Gr} variations of Hodge structure. The second property
is a theorem of Tian~\cite{T1} (See also Todorov~\cite{To}). The third one
is the compactification theorem of Viehweg~\cite{V2} and the forth
property can be verified by  the Nilpotent Orbit theorem of
Schmid~\cite{Schmid}.
\end{rem}

When we study the boundary behavior of the Weil-Petersson metric on the
moduli spaces, we need to analyze the potential 
$(\Omega,\bar\Omega)$, where $(\,,\,)$ is the polarization
of the Hodge structures.
We recall here the result of the  potential of Weil-Petersson
metric on one dimensional slice of moduli spaces.

Let $\Delta^*$ be a one dimensional parameter space of
a family of polarized Calabi-Yau manifolds. Let $\Omega$
be a section of the first Hodge bundle $F^n$. Then by
the Nilpotent Orbit theorem of Schmid \cite{Schmid},
after
a possible base change,  we have
\[
\Omega=e^{\bb N\log\frac 1z}A(z),
\]
where $N$ is the nilpotent operator, $N^{n+1}=0$ for
$n$ being the dimension of the Calabi-Yau manifolds, 
and 
\begin{equation}\label{Aseries}
A(z)=A_0+A_1z+\cdots
\end{equation}
is a
vector-valued convergent power series with the
convergent radius $\delta>0$.
 Let  
\[
f_{k,l}(z)=z^k(\log\frac 1z)^l,
\]
for any $k,l\geq 0$. Then we can write $\Omega$ as the
convergent series
\begin{equation}\label{bseries}
\Omega=\sum_{k,l}A_{k,l} z^k(\log\frac 1z)^l=\sum_{k,l}
A_{k,l}f_{k,l}.
\end{equation}

Define $\deg f_{k,l}=k-\frac{l}{n+1}$. Then
we have the following lemma (Lemma 7.1 of \cite{LS-1}):

\begin{lemma}\label{sp1}
The convergence of ~\eqref{bseries} is in the
$C^\infty$ sense. Furthermore, we have
\begin{equation}\label{esti}
||\Omega-\sum_{\deg
f_{k,l}\leq\mu}A_{k,l}f_{k,l}||_{C^s}\leq C
r^{k_0-s}(\log\frac 1r)^{l_0},
\end{equation}
where $r=|z|$, $k_0,l_0$ are the unique pair of nonnegative
integers such that
$l_0\leq n$,
$k_0-\frac{l_0}{n+1}>\mu$  and for any pair of
integers $k',l'$ with 
$k'-\frac{l'}{n+1}>\mu$ we have
$k'-\frac{l'}{n+1}\geq k_0-\frac{l_0}{n+1}$. 
$C$ is  a constant depending only on $k_0, l_0, \mu$ and
$\Omega$.
\end{lemma}

In order to estimate the volume form and the Chern classes of the 
Weil-Petersson metric, we also need the following Strominger's 
formula (Theorem 3.1 of \cite{LS-1}). 
\begin{theorem}
Let $(g_{i\bar j})_{m \times m}$ be the Weil-Petersson metric and let 
$D_jD_i\Omega$ be the projection of $\partial_j\partial_i\Omega$ onto 
$H^{n-2,2}$. Then the curvature tensor of the Weil-Petersson metric is 
given by
\[
R_{i\bar j k\bar l}=g_{i\bar j}g_{k\bar l}+g_{i\bar l}g_{k\bar j}
-\frac{(D_kD_i\Omega,\bar{D_lD_j\Omega)}}{(\Omega,\bar\Omega)}
\]
for $1 \leq i,j,k,l \leq m$. 
\end{theorem}

In order to bound the Weil-Petersson metric and its curvature,
we  need the Hodge metric, which was defined
in~\cite{Lu3}.
\begin{theorem}
Let $D$ be the classifying space. 
The invariant Hermitian metric of $D$ restricts to $ M$ is a 
K\"ahler
metric called the Hodge metric. Let $\omega_H$ be its K\"ahler form. 
Then we have
\begin{enumerate}
\item The bisectional curvature of $\omega_H$ is nonpositive;
\item $Ric(\omega_H)<\alpha\omega_H<0$ for some negative constant 
$\alpha$ which only depends on the dimension;
\item The holomorphic sectional curvature of $\omega_H$ is bounded
above by $\alpha$;
\item $2\omega_{WP} \leq \omega_H$ and $-\omega_H \leq
Ric(\omega_{WP})\leq \omega_H$.
\end{enumerate}
\end{theorem}
The Hodge metric is useful because it gives us the convexity of
the volume form of the Weil-Petersson metric. Furthermore,
using the Hodge metric together with the Schwarz-Yau Lemma, we can
control the volume of the Weil-Petersson metric and the Hodge
metric.
\begin{theorem}
Let $M$ be the moduli space of polarized Calabi-Yau $n$-folds.
Then the volume of 
Weil-Petersson subvariety $M_1$ of $M$ equipped with
the Weil-Petersson metric or the Hodge metric is finite. 
\end{theorem}

\section{Cut-off functions}\label{cut}
The main result of this section is to prove the
existence of the ``good'' cut-off functions so that our
estimates can go through.

We assume that $\bar M$ is an $m$-dimensional compact \ka manifold of
dimension $m$ and $D$ be a divisor of $\bar M$ with
normal crossing so that $M=\bar M\backslash D$. 
We are going to prove that there is a complete \ka metric
on
$M$ such that it is asymptotical to  the Poincar\'e
metric near infinity. We call 
this metric the global Poincar\'e metric or simply
the Poincar\'e metric. We use $\omega_P$ to denote
its \ka form.

The result is well known, for example, in~\cite{JY1}.
For the sake of completeness and for the setting of
notations, we sketch the proof in Lemma~\ref{lem31}.

We let $\bar M=U_1\cup\cdots \cup U_t\cdots \cup U_s$
be covered by local coordinate charts where $1\leq t<s$.
Without loss  generality, we assume that $(\bar
U_{t+1}\cup\cdots
\cup \bar U_s)\cap D=\emptyset$. For each $1\leq\alpha\leq
t$, we assume that there is an $l_\alpha$ such that each $U_\alpha\backslash D
=(\Delta^{\ast})^{l_{\alpha}}\times\Delta^{m-l_{\alpha}}$ and 
on each $U_\alpha$,
$D$ is defined by 
\[
z_1^\alpha\cdots z_{l_\alpha}^\alpha=0.
\]
Let $\{\psi_\alpha\}_{1\leq\alpha\leq s}$ be the 
partition of the unity subordinated to the cover
$\{U_\alpha\}_{1\leq\alpha\leq s}$. Let $\omega$ be a
\ka  metric of $\bar M$ and let
$C$ be a large constant. Define 
\[
\omega_P=C\omega-\sum_{\alpha=1}^t\sqrt{-1}\pa\bar\pa
\left(\psi_\alpha\sum_{j=1}^{l_\alpha}\log\log\frac{1}
{|z_j^\alpha|^2}\right).
\]
Then we have

\begin{lemma}\label{lem31}
For $C$ large enough, $\omega_P$ defines a complete
 metric on $M$ with finite volume and bounded
curvature. Furthermore, there is a constant $C_1$ such
that
\[
\frac{1}{C_1}\omega_0^\alpha\leq\omega_P\leq C_1
\omega_0^\alpha
\]
for any $1\leq\alpha \leq t$, where
$\omega^\alpha_0$ is the local Poincar\'e metric, defined
by
\[
\omega^\alpha_0=\sum_{i=1}^{l_{\alpha}}\sqrt{-1}\frac{1}{r_i^2(\log\frac{1}
{r_i})^2}dz^\alpha_i\wedge d\bar z^\alpha_i
+\sum_{i=l_{\alpha}+1}^m\sqrt{-1}dz_i^\alpha\wedge 
d\bar z_i^\alpha.
\]
\end{lemma}

{\bf Proof.} This follows from a straightforward
computation.

\qed

The main result of this section is the following:

\begin{theorem}\label{thm41}
Let $\eps>0$. 
Then there is a function  $\rho_\eps$
such that
\begin{enumerate}
\item
$0\leq\rho_\eps\leq 1$;
\item
For any open neighborhood $V$ of $D$ in $\bar M$, there is
$\eps>0$ such that ${\rm supp}(1-\rho_\eps)\subset V$;
\item For each $\eps>0$, there is a neighborhood
$V_1$ of $D$ such that $\rho_\eps|_{V_1}\equiv 0$;
\item $\rho_{\eps'}\geq\rho_\eps$ for $\eps'\leq \eps$;
\item There is a constant $C$, independent of $\eps$ such
that
\[
-C\omega_P\leq\sqrt{-1}\pa\bar\pa\rho_\eps\leq C\omega_P.
\]
\end{enumerate}
\end{theorem}

{\bf Proof.} The key observation can be explained as
follows. For the  unit ball in $\R^n$, if
we construct a smooth function which is $0$ on $B_\eps(0)$
and $1$ outside $B_{2\eps}(0)$, then the second derivative
of the function is in general of order $1/\eps^2$.
However, in the two dimensional case, if we consider the
Laplacian of the function, then it is possible to lower
the order of the second derivative.

Define a decreasing smooth function $\phi: \R\rightarrow \R,
0\leq\phi\leq 1$ as follows:
\begin{equation}\label{eqn41}
\phi(t)=\left\{
\begin{array}{ll}
0& t\geq 1;\\
1& t\leq 0.
\end{array}
\right.
\end{equation}
and assume that $|\phi'|+|\phi''|<10$.
Define a function $\phi_\eps$ on the unit ball of the
complex plane as 
\[
\phi_\eps(z)=\phi\left(\frac{(\log\frac
1r)^{-1}-\eps}{\eps}\right),
\]
where $r$ is the Euclidean norm of the complex variable $z$.
Then we have the following
\begin{align*}
&\pa_z\phi_\eps=\frac{1}{2\eps}\phi'\frac{1}
{z(\log\frac 1r)^2},\\
&\pa_z\bar\pa_z\phi_\eps=\frac{1}{4\eps^2}\phi''\frac{1}{r^2(\log\frac
1r)^4}+\frac{1}{2\eps}\phi'\frac{1}{r^2(\log\frac 1r)^3}.
\end{align*}
Thus we have

\begin{lemma}\label{lem42}
Using the same notation as above, we have
\begin{align*}
&|\pa\phi_\eps|\leq \frac{10}{r(\log\frac 1r)}\\
-10\sqrt{-1}\frac{1}{r^2(\log\frac 1r)^2}dz\wedge d\bar z
&\leq \sqrt{-1}\pa\bar\pa\phi_\eps
\leq 10\sqrt{-1}\frac{1}{r^2(\log\frac 1r)^2}dz\wedge
d\bar z,
\end{align*}
where the norm is with respect to the Euclidean 
metric on $\C$.
Furthermore, we have $${\rm supp} (\pa\phi_\eps)\subset
B_{e^{-\frac{1}{2\epsilon}}}-B_{e^{-\frac{1}{\epsilon}}}$$.
\end{lemma}

\qed

For $1\leq\alpha\leq t$ and $\eps>0$ small enough, 
let 
\[
\phi_\eps^\alpha(z_1^\alpha,\cdots,z_m^\alpha)
=\Pi_{i=1}^{l_\alpha}(1-\phi_\eps(z_i^\alpha))
\]
on $U_\alpha$.
Then using the Lemma~\ref{lem42}, we have
\begin{align}\label{eqn43}
\begin{split}
& |\pa\phi_\eps^\alpha|\leq C\sum_{i=1}^{l_\alpha}
\left |\frac{1}{r_i(\log\frac{1}{r_i})}\right |,\\
& -C\omega^\alpha_0\leq\sqrt{-1}\pa\bar\pa
\phi_\eps^\alpha\leq C\omega_0^\alpha
\end{split}
\end{align}
for some constant $C$.

We define
\[
\rho_\eps=1-\sum_{\alpha}^t\psi_{\alpha}\phi_{\eps}^{\alpha}.
\]
We can verify that $\rho_\eps$ satisfies all the properties
in Theorem~\ref{thm41} by a straightforward computation.

\qed

We finish this section by stating the Schwarz-Yau Lemma~\cite{Y3}
in our context.  The result will be used repeatedly in the rest of this paper:

\begin{prop}\label{SY}
Let $M$ be a smooth Weil-Petersson variety whose compactification $\bar M$ is also smooth. Let $D=\bar M-M$ be a divisor of normal crossings. Let $\omega$ is a K\"ahler
metric of $M$ such that its holomorphic sectional curvature is less than a constant $-\delta$ for $\delta>0$. Then there is a constant $C$ such that
\[
\omega\leq C\omega_P,
\]
where $\omega_P$ is the global Poincar\'e metric.
\end{prop}

\section{Extension of the Hodge bundles}\label{ext}
Let $M$ be a Weil-Petersson variety and $\bar M$ be its
smooth compactification. We assume that $Y=\bar 
M-M$ be the divisor of normal crossings. 
Then it is known from~\cite[page 235]{Schmid} that the
Hodge bundles
$F^n,\cdots, F^0$ can be extended to  coherent sheaves
over $\bar M$. Furthermore, if we assume that  every
element of the monodromy group is unipotent,
then the coherent sheaves are in fact vector bundles
over $\bar M$.

The particular extension of the bundles over $\bar M$
is defined by the Nilpotent Orbit theorem of
Schmid. In \S 5 and \S 6, we  show that such an
extension is the one that we can control the growth
of the Weil-Petersson metric  and its curvatures at
infinity. ~\footnote{In fact, in the case of the
moduli space of Riemann surfaces, there is only one
extension of the bundles that is ``good'' in the
sense of Mumford~\cite[Section 1]
{Mumford1}(cf.~\cite{Tra1}). In our case, even if
the  extension may  not be ``good'', it is the best
possible extension we can get.} Thus although we
know the extension exists by the Nilpotent Orbit
theorem, we must write out explicitly the local
transition functions.

The first result of this section is the following lemma
which is essentially due to Kawamata~\cite{kawa3}.
We formulate it in the language of  Weil-Petersson
geometry. 

\begin{lemma}\label{lem41}
Let $M$ be a Weil-Petersson variety with $\bar M$ being its
compactification. We assume that $\bar M$ is smooth, projective and
$Y=\bar M-M$ is a divisor of normal crossings. Then there is 
a divisor $Y_1$ of $M$ of normal crossings such that
there is a finite covering $M_1$ of the variety $M-Y_1$
with the following properties:
\begin{enumerate}
\item $M_1$ is a Weil-Petersson manifold;
\item The elements of the monodromy group of $M_1$
are unipotent.
\end{enumerate}
\end{lemma}

{\bf Proof.} We first observe that if we remove
any divisor $Y_1$ from $M$, then $M-Y_1$ is still a 
Weil-Petersson variety. For this reason, without loss
generality, we can assume that the  Weil-Petersson variety 
$M$ is
actually a Weil-Petersson  manifold.

Now we use the idea of Kawamata~\cite{kawa3}. Let $T$ be a
monodromy operator which is not unipotent. Let
$T=\gamma_{s}\gamma_{u}$ be the decomposition of $T$ into its
semi-simple part and its unipotent part. By the theorem of Borel, 
there is an integer $m$ such that $\gamma_s^m=1$. Let $L$ be an
ample line bundle of $\bar M$. 
Let $Y+Y_1=\sum D_j$ be the decomposition of the divisor
$Y+Y_1$ into irreducible pieces,
where $Y_1$ is the divisor containing the singular locus of $M$.
 We assume that the monodromy
operator $T$ is generated by $U\backslash D_1$,
where $U$ is a  neighborhood of $D_1$.
Assume that $s$ is large enough
such that the bundle $L^s(-D_1)$ is very ample. By taking the
$m$-th root of the sections of $L^s(-D_1)$ we get a variety
$M_1$ such that outside a possible divisor, it is a finite
covering space of $M$. $M_1$ may have some singularities. However,
we can always remove those 
divisors containing singularities to get a smooth manifold. 
The explicit construction of $M_1$ is as follows: let $S_0$
be a generic section of the line bundle $L^s(-D_1)$. $S_0$
is generic in the sense that on $D_i$ with $i\neq 1$, $S_0$
is not identically zero and $dS_0\neq 0$ generically on $D_1$. We
can extend
$S_0$ to a basis
$S_0,\cdots, S_t$ of $H^0(L^s(-D_1))$ such that the basis
defines an embedding
\[
\sigma: 
\bar M\rightarrow {\mathbb C}P^t,\quad
x\mapsto [S_0,S_1,\cdots, S_t].
\]
We consider the map $\pi: CP^t\rightarrow CP^t$ by
$[Z_0,\cdots,Z_t]\mapsto [Z_0^m,Z_1,\cdots,Z_t]$.
It is a holomorphic $m$-branched covering map. Let $Z$ be the
pre-image of $\bar M$ under $\pi$. Then $Z$ is a projective
variety. Let
$Z_{\rm reg}$ be the smooth points of $Z$.  Define
$M'=Z_{\rm reg}\cap \{Z_0=0\}$ and let $\tilde Z$ be the
desingularization of $Z$ along the divisors $Z\backslash M'$.

The Hodge
bundles
can be pulled back to the  manifold $M'$. At any neighborhood
$(\Delta^*)^l\times \Delta^{k-l}$ of $\bar M\backslash M$, the
transform of
$(M,\bar M)$ to
$(M', \tilde Z)$ is the $m$-branched covering defined by
$z_1\mapsto
\sqrt[m]{z_1}$, where $z_1=0$ is corresponding to the divisor
$D_1$.  Evidently, the monodromy operator $T$ is transformed to
$T^m$, which becomes  a unipotent operator. 

One can observe that if $T'$ is a unipotent  monodromy operator,
then
 under the transform $(M, \bar M)\mapsto (M',\tilde Z)$, $T'$ is
still unipotent.   
Since there are only finitely many irreducible components
of $D$, there are only finitely many monodromy operators which are not unipotent.
Thus by finitely many  transforms, we can
get a Weil-Petersson  manifold $M_1$ on which all
monodromy operators are unipotent.

\qed

\begin{rem} In general, a Weil-Petersson variety $M$ may have orbifold singularities.
However, we let $M'$ be the regular part of $M$. Then $M'$ is a manifold that has 
a smooth compacification $\bar M'$. By Lemma ~\ref{lem41}, up to a finite covering, we can assume that
the elements of the monodromy group are unipotent. Let $M''\rightarrow M'$ be the finite covering defined by Lemma~\ref{lem41}, then there is an integer $s$ such that
\[
{\rm vol} (M'')=s\,{\rm vol} (M')
\]
and
\[
\int_{M''}c_1(\omega_{WP})^k\wedge\omega_{WP}^l
=s
\int_{M'}c_1(\omega_{WP})^k\wedge\omega_{WP}^l.
\]
Thus
From now on, we will prove our results under the additional assumptions that  all monodromy operators
are unipotent, $M$, $\bar M$ are smooth, and the divisor  $Y=\bar{ M}\setminus
{M}$ is  of 
normal crossings.
\end{rem}

We  write out explicitly
the extension of the Hodge bundles in terms of the local 
coordinates.

 A bundle $F^p$ over ${M}$ is
equivalent to an open cover $\{V_\alpha\}$ of $M$ with
transition functions $g_{\alpha\beta}: V_\alpha \cap V_\beta \to
GL(s,\mathbb C)$ with $s=
\text{rank}\ F^p$. We  assume that $\{V_\alpha\}$
is a countable, locally finite cover. Let $\phi:
 M \to \Gamma\backslash D$ be the period map,
where $\Gamma$ is the monodromy group.
Since $Y \subset \bar{ M}$ is compact, we can
take a finite cover $\{ U_\alpha \}_{\alpha=1}^{t}$ of $Y$ in
$\bar{M}$ such that
\begin{enumerate}
\item Each $U_\alpha$ is bi-holomorphic to
$\Delta^m$, the polydisc on $\C^m$.
\item $\bigcup_1^t U_\alpha$ contains $\bar S$,
where $S$ is a
neighbhood  of $Y$ in $\bar{ M}$.
\item On each $U_\alpha$ with local coordinates
$z_1^\alpha,\cdots,z_m^\alpha$, the divisor $Y \cap
U_\alpha$ is given by $z_1^\alpha\cdots
z_{l_\alpha}^\alpha=0$ for some $l_{\alpha} \in \{
1,\cdots,m\}$.
\end{enumerate}

Since $\cup\{ U_\alpha \} \cup \{  V_\alpha \}$ is an
open cover of $\bar{M}$, we can take a finite
subcover. After refinement, we can take
$V_1,\cdots,V_r \in \{ V_\alpha \}$ such that
\begin{enumerate}
\item $\bigcup_1^t U_\alpha \cup \bigcup_1^r
V_\alpha=\bar{M}$.
\item $\bar{V_\alpha} \cap Y=\emptyset$ for $1\leq\alpha \leq r$.
\end{enumerate}

Now, for each $U_\alpha=\Delta^m$, let
$U_\alpha^\ast=U_\alpha\setminus
Y=(\Delta^\ast)^{l_\alpha}\times \Delta^{m-l_\alpha}$.
We cut $U_\alpha^\ast$ into open conical parts such that
each conical part is a product of  discs and open sectors 
with small angles. We 
write $U_\alpha^\ast=\bigcup U_{\alpha,i}$ where each
$U_{\alpha,i}$ is a conical domain.
~\footnote{For example, we can define each
$U_{\alpha,i}$ to be
$U_{\alpha,i}=\{(z_1^\alpha,\cdots,z_m^\alpha)\in U_\alpha^*|\,
{\rm arg} z_i^\alpha\in (a_i,b_i), 1\leq i\leq
l_\alpha\}$,
where $|b_i-a_i|$ is small.}

For each $U_{\alpha}^\ast$, the universal covering space
of
$\widetilde{U}_\alpha^\ast$ is $U^{l}\times \Delta^{m-l}$ where
$U$ is the upper half-plane and $l=l_\alpha$. The
natural projection $p: \widetilde{U}_\alpha^\ast\to
U_{\alpha}^\ast$ is given by
$p(w_1,\cdots,w_l,\cdots,w_m)=(e^{2\pi
iw_1},\cdots,e^{2\pi iw_l},w_{l+1},\cdots,w_m)$.
Let $\tilde M$ be the universal covering space of $M$.
Then we always have the lifting $\tilde\phi:\tilde M\rightarrow D$.
Locally, it is given by 
 the map
$\widetilde{\phi}_\alpha:\widetilde{U}_\alpha^\ast \to D$
corresponding to 
$\phi_\alpha:U_{\alpha}^\ast \to \Gamma\backslash
D$, and we have  the following commutative diagram
\[
\begin{CD}
\widetilde{U}_\alpha^\ast @> >>D\\
@VVV        @VVV\\
U_\alpha@ >>> \Gamma\backslash D
\end{CD}.
\]
 Furthermore, $\tilde\phi_\alpha$ and
$\tilde\phi_\beta$ are compatible if $U_\alpha\cap
U_\beta\neq\emptyset$.

By the Nilpotent Orbit theorem of Schmid, for each
$j=1,\cdots,l$, there is a monodromy transform
$T_j$ such that
\[
\widetilde{\phi}_\alpha(w_1,\cdots,w_j+1,\cdots,w_m)=T_j
\circ
\widetilde{\phi}_\alpha(w_1,\cdots,w_m).
\]
We have $T_jT_k=T_kT_j$
for $1 \leq j,k\leq l$. By Lemma~\ref{lem41}, we
assume that the semisimple part of each $T_j$
is $1$. Let
\[
\widetilde{\psi}_\alpha(w_1,\cdots,w_m)
=e^{-\sum_{1}^{l}w_jN_j}
\widetilde\phi_\alpha(w_1,\cdots,w_m),
\]
where $N_j=\log T_j$ are the nilpotent operators. Clearly
$\widetilde{\psi}_\alpha$ is invariant under the transform
$w_j \mapsto w_j+1$ for $j=1,\cdots,l$. So
$\widetilde{\psi}_\alpha$ descends to a map $\psi_\alpha$ from
$U_\alpha^\ast$ to the complex dual $\check D$ of
$D$. By the Nilpotent Orbit theorem, $\psi_\alpha$
can be holomorphically extended to
$U_\alpha$.

Since the neighbhood $S$ of $Y$ in $\bar{M}$ satisfies
the condition that $\bar S \cap V_{\beta}=\emptyset$ for each
$V_\beta$, and for each $U_\alpha$, there is a $\sigma_\alpha>0$
such that $\bigcup_{j=1}^{l_\alpha} \Delta_{1}^\ast\times
\cdots\times\Delta_j^\ast(\sigma_\alpha)\times\cdots \times
\Delta^*_{l_\alpha}\times \Delta^{m-l_\alpha} \subset S\cap
U_\alpha$, we know that $(\Delta_{1}^\ast\times
\cdots\times\Delta_j^\ast(\sigma_\alpha)\times\cdots \times
\Delta_{l_\alpha}^*\times \Delta^{m-l_\alpha})\cap V_\beta=\emptyset$
for any $1 \leq j \leq l_\alpha$ and $V_\beta$. Let
$U_{\alpha,0}=(\Delta^\ast(\sigma_\alpha))^{l_\alpha}
\times\Delta^{m-l_\alpha}$. 
If  ${\rm diam\,}U_\alpha\rightarrow 0$,
then  we can assume that $Y \subset \bigcup
U_{\alpha,0}$ and $U_{\alpha,0} \cap
V_\beta=\emptyset$.

Now, $\{ U_{\alpha,0} \} \cup \{ U_{\alpha,i} \} \cup \{ V_{\beta}
\}$ give an open cover of $\bar{ M}$. Clearly, on each
$U_{\alpha,i}$ and each $V_\beta$, the Hodge bundles are
trivial. For each $0 \leq p \leq n$, let
$s=\text{rank}\ F^p$. We will extend $F^p$ to
$\bar{ M}$ using the cover $\{ U_{\alpha,0} \}
\cup \{ U_{\alpha,i} \} \cup \{ V_{\beta} \}$ of
$\bar{M}$. Since $F^p$ over each
simply connected set of $U^\ast_\alpha$ are trivialized by
$e^{\frac{\sqrt{-1}}{2\pi}\sum_{1}^{l}
\log\frac{1}{z_j}N_j}\psi_\alpha$, we know that
$g_{_{U_{\alpha,i},U_{\alpha,j}}}=I$ where $I=I_s$ is the identity
matrix of rank $s$. Since $F^p$ is trivial over $U_{\alpha,i}$ and
$V_\beta$, we know that the transition functions 
$g_{_{U_{\alpha,i},V_{\beta}}}$ and
$g_{ij}^{\alpha\beta}=g_{_{U_{\alpha,i},U_{\beta,j}}}$
are given. 
Since $U_{\alpha,0}\cap V_\beta=\emptyset$, 
to extend $F^p$ to $\bar M$, we only need to 
define $g_{_{U_{\alpha,0},U_{\alpha,j}}}$,
$g_{_{U_{\alpha,0},U_{\beta,j}}}$, and 
$g^{\alpha\beta}=g_{_{U_{\alpha,0},U_{\beta,0}}}$ . 

We define
$g_{_{U_{\alpha,0},U_{\alpha,j}}}=I$. Since $U_{\alpha,0}\setminus Y 
\subset \bigcup U_{\alpha,i}$, we know that for each $q \in U_{\alpha,0}
\cap U_{\beta,j}$, there is a $U_{\alpha,i}$ such that $q \in
U_{\alpha,i}$. So naturally we define 
$g_{_{U_{\alpha,0},U_{\beta,j}}}$ 
on $U_{\alpha,0}\cap U_{\beta,j}\cap U_{\alpha,i}$
to be the restriction of 
$g_{_{U_{\alpha,i},U_{\beta,j}}}$ on $U_{\alpha,0}
\cap U_{\beta,j}\cap U_{\alpha,i}$.

Let $U_{\alpha,i'}\cap U_{\alpha,i}\cap U_{\beta,j}
\cap U_{\alpha,0}\neq\emptyset$. Then since
$g_{_{U_{\alpha,i},U_{\alpha,i'}}}=I$, we have
$g_{_{U_{\alpha,i},U_{\beta,j}}}
=g_{_{U_{\alpha,i'},U_{\beta,j}}}$ on
$U_{\alpha,i'}\cap U_{\alpha,i}\cap U_{\beta,j}
\cap U_{\alpha,0}$. Thus $g_{_{U_{\alpha,0},
U_{\beta, j}}}$ are well-defined.

It is more difficult  to define $g^{\alpha\beta}$.
For each 
$q\in U_{\alpha,0}\cap U_{\beta,0}$, if $p \notin Y$, then there exist
$i$ and $j$ such that $q\in U_{\alpha,i}\cap U_{\beta,j}$. We define 
$g^{\alpha\beta}(q)=g_{ij}^{\alpha\beta}(q)
=g_{_{U_{\alpha,i},U_{\beta,j}}}$. As long as these
transition functions are well-defined, it is trivial to check the
compatibility conditions. Thus we need to prove

\begin{claim}\label{claim1} Using the above
notations, we have
\begin{enumerate}
\item $g^{\alpha\beta}$ are well-defined
on $U_{\alpha,0}\cap U_{\beta,0}\backslash Y$.
\item $g^{\alpha\beta}$ can be extended to  $
U_{\alpha,0}\cap U_{\beta,0}$.
\end{enumerate}
\end{claim}
\begin{rem}\label{bdddd}
Fixing a pair of charts $U_{\alpha,0}$ and
$U_{\beta,0}$ such that their  intersection is
non-empty, we can choose compatible
coordinates on $U_{\alpha,0}$ and $U_{\beta,0}$.
Assume
$Y_{\alpha\beta}=U_{\alpha,0}\cap U_{\beta,0}\cap
Y$ is of codimension $r$ and is given by
$z_{1}^\alpha
\cdots z_r^\alpha=0$ and 
$z_{1}^\beta \cdots z_r^\beta=0$ on
$U_{\alpha,0}$ and $U_{\beta,0}$, respectively.
Since these two equations define the same variety
$Y_{\alpha\beta}$, we can identify the variables
pairwisely in the following way:  for each
$i=1,\cdots,r$,  choose a point $q \in
Y_{\alpha\beta}$ such that $z_i^\alpha(q)=0$ and
$z_j^\alpha(q)\ne 0$ for $j\ne i$. Clearly there is a $t$ with $1 \leq t
\leq r$ such that $z_t^\beta(q)=0$ and $z_u^\beta(q)\ne 0$ for $u\ne t$. 
Without loss of generality, we assume that $t=i$.
Thus
$\frac{z_i^\beta}{z_i^\alpha}$ is non-zero for $1\leq i \leq r$.
Furthermore, by slightly shrinking $U_{\alpha,0}$ and $U_{\beta,0}$ we
can assume that $z_j^\alpha$ and $z_j^\beta$ are
bounded and are bounded away from $0$ on
$U_{\alpha,0}\cap U_{\beta,0}$ for $r+1 \leq j \leq
l_\alpha$.
\end{rem} 

Let $D$ be the classifying space. By definition,
this means that $D$ is the space of decompositions
of a fixed vector space $H$ satisfying
the Hodge-Riemann relations. 
To prove the first assertion of Claim~\ref{claim1},
we fix a
$q\in (U_{\alpha,0}\cap U_{\beta,0})
\setminus Y$. Assume that $q\in U_{\alpha,i}\cap U_{\beta,j}$. 
Let $v_1,\cdots,
v_b$ be a fixed  basis of $H$.
Using this basis, we can identify a map $\rho$ to
the classifying space (or to its compact dual) with
 a sequence of matrices.
Let $\pi_p\circ
\rho\in F^p$.
Then $\pi_p\circ\rho$ can be represented by
a $b\times s$ matrix-valued function,
where $s=\dim\, F^p$ as above.

Since $
e^{\frac{\sqrt{-1}}{2\pi}\sum_{1}^{l_\alpha}
\log\frac{1}{z^\alpha_i}N_i}\psi_\alpha$ gives
local trivialization of
$F^p$ over $U_{\alpha,i}$ and $
e^{\frac{\sqrt{-1}}{2\pi}\sum_{1}^{l_\beta}
\log\frac{1}{z^\beta_j}N_j}\psi_\beta$ gives local
trivialization of
$F^p$ over $U_{\beta,j}$, 
 both of them can be represented by  $b\times s$
matrices.
Fix a lift  and lift $\phi_\alpha$ and
$\phi_\beta$ to the universal covering spaces 
$\widetilde{U}_\alpha^\ast$ and $\widetilde{U}_\beta^\ast$ respectively.
We have that, as matrices 
\begin{eqnarray*}
\widetilde{\phi}_\alpha(w_1^\alpha,\cdots,w_m^\alpha)=
\widetilde{\phi}_\beta(w_1^\beta,\cdots,w_m^\beta)g_{ij}^{\alpha\beta}
(w_1^\beta,\cdots,w_m^\beta).
\end{eqnarray*}
For each $1 \leq i \leq r$, let $T_i$ be the monodromy transformation
such that 
\begin{eqnarray*}
\widetilde{\phi}_\alpha(w_1^\alpha,\cdots,w_i^\alpha+1,\cdots,w_m^\alpha)
=T_i\, 
\widetilde{\phi}_\alpha(w_1^\alpha,\cdots,w_i^\alpha,\cdots,w_m^\alpha)
\end{eqnarray*}
and 
\begin{eqnarray*}
\widetilde{\phi}_\beta(w_1^\beta,\cdots,w_i^\beta+1,\cdots,w_m^\beta)
=T_i\,
\widetilde{\phi}_\beta(w_1^\beta,\cdots,w_i^\beta,\cdots,w_m^\beta).
\end{eqnarray*}
This is true  since $T_i$ is corresponding to the
same simple loop of 
$(U_{\alpha,0}\cap U_{\beta,0}) \setminus Y$.
Combining the above three formulae, we have 
\begin{align*}
\begin{split}
&T_i\circ 
\widetilde{\phi}_\beta(w_1^\beta,\cdots,w_i^\beta,\cdots,w_m^\beta)
g_{ij}^{\alpha\beta}(w_1^\beta,\cdots,w_i^\beta+1,\cdots,w_m^\beta)\\
=& \widetilde{\phi}_\beta(w_1^\beta,\cdots,w_i^\beta+1,\cdots,w_m^\beta)
g_{ij}^{\alpha\beta}(w_1^\beta,\cdots,w_i^\beta+1,\cdots,w_m^\beta)\\
=&
\widetilde{\phi}_\alpha(w_1^\alpha,\cdots,w_i^\alpha+1,\cdots,w_m^\alpha)
\\
=&T_i\circ 
\widetilde{\phi}_\alpha(w_1^\alpha,\cdots,w_i^\alpha,\cdots,w_m^\alpha)\\
=&T_i\circ 
\widetilde{\phi}_\beta(w_1^\beta,\cdots,w_i^\beta,\cdots,w_m^\beta)
g_{ij}^{\alpha\beta}(w_1^\beta,\cdots,w_i^\beta,
\cdots,w_m^\beta),
\end{split}
\end{align*}
which implies
\begin{align*}
\begin{split}
g_{ij}^{\alpha\beta}(w_1^\beta,
\cdots,w_i^\beta+1,\cdots,w_m^\beta)
=
g_{ij}^{\alpha\beta}(w_1^\beta,\cdots,w_i^\beta,\cdots,w_m^\beta).
\end{split}
\end{align*}
 So 
$g_{ij}^{\alpha\beta}$ is invariant under the deck 
transformations of the
universal covering of $U_{\alpha,0}\cap U_{\beta,0}
\setminus Y$ which
implies that it descends to a function on 
$U_{\alpha,0}\cap U_{\beta,0}\setminus Y$. This proved that 
$g_{ij}^{\alpha\beta}$ is well-defined. 

Now we prove the second assertion of
Claim~\ref{claim1}. As stated in the Remark
\ref{bdddd}, if we let
$\xi_i=\frac{z_i^\beta}{z_i^\alpha}$, we know that
$\xi_i$ is bounded and bounded away from $0$ for
$1\leq i \leq r$. Since 
\begin{eqnarray}\label{1111}
e^{\frac{\sqrt{-1}}{2\pi}\sum_{1}^{l_\alpha}
\log\frac{1}{z^\alpha_i}N_i}\psi_\alpha
=e^{\frac{\sqrt{-1}}{2\pi}\sum_{1}^{l_\beta}
\log\frac{1}{z^\beta_i}N_i}\psi_\beta g_{ij}^{\alpha\beta}
\end{eqnarray} 
on $U_{\alpha,0}\cap U_{\beta,0}$, we have 
\begin{eqnarray}\label{2222}
e^{\frac{\sqrt{-1}}{2\pi}\sum_{r+1}^{l_\alpha}
\log\frac{1}{z^\alpha_i}N_i}\psi_\alpha
=e^{\frac{\sqrt{-1}}{2\pi}(\sum_1^r\log\frac{1}{\xi_j}N_j+
\sum_{r+1}^{l_\beta}
\log\frac{1}{z^\beta_i}N_i)}\psi_\beta g_{ij}^{\alpha\beta}.
\end{eqnarray}
By the definition of $\xi_j$ and by Remark \ref{bdddd} we know that 
$\log\frac{1}{z^\alpha_i}$ for $r+1 \leq i \leq l_\alpha$,
$\log\frac{1}{z^\beta_i}$ for $r+1 \leq i \leq l_\beta$ and
$\log\frac{1}{\xi_j}$ for $1 \leq i \leq r$ are bounded. By the Nilpotent
Orbit theorem we know that 
$\psi_\alpha$ and $\psi_\beta$ can be holomorphically extended to
$U_{\alpha,0}$ and $U_{\beta,0}$ respectively and image of $\psi_\alpha$
and $\psi_\beta$ restricted to $Y$ lie in the complex dual of $D$. So we
can find a non-singular $r \times r$ minors $A_\alpha$ of
$\psi_\alpha$. Denote the corresponding minor of $\psi_\beta$ by
$A_\beta$. By
\eqref{2222}, since $N_i$ are fixed nilpotent operators, we know that 
$g_{ij}^{\alpha\beta}=C A_{\alpha} A_{\beta}^{-1}$ where 
$C=e^{\frac{\sqrt{-1}}{2\pi}\sum_{r+1}^{l_\alpha}
\log\frac{1}{z^\alpha_i}N_i
-\frac{\sqrt{-1}}{2\pi}\sum_1^r\log\frac{1}{\xi_j}N_j
-\frac{\sqrt{-1}}{2\pi}\sum_{r+1}^{l_\beta}
\log\frac{1}{z^\beta_i}N_i}$ is inveritible. Since $A_\beta$
is also bounded, we know that $g_{ij}^{\alpha\beta}$ is bounded away
from $0$. The same argument works for $(g_{ij}^{\alpha\beta})^{-1}$. So 
$g_{ij}^{\alpha\beta}$ is also bounded. This implies that
$g_{ij}^{\alpha\beta}$ can be extended 
to $U_{\alpha,0}\cap U_{\beta,0}$ which  finishes the
construction of the extensions of the Hodge bundles. 

\begin{rem}
The tangent bundle of $ M$ can be identified with 
$F^{n-1}/ F^n$. Since all the Hodge bundles can be extended in the
above canonical  way, the tangent  bundle is also canonically
extended.
\end{rem}

\section{Volume of the moduli space}\label{vol}
In this section we prove that the volume of the Weil-Petersson metric is 
a rational number. Let $\Delta_r$ be the disk 
of radius $r$ and $\Delta^{\ast}_{r}$ be the 
punctured disk of radius $r$ in $\mathbb{C}$.  
Let $V^{k}_{r}=(\Delta^{\ast}_{r})^k\times
\Delta_{r}^{m-k}$ for $1 \leq k \leq m$. 

Assume that a chart of ${M}$ near the boundary is $V_{1}^{k}$ and 
the Weil-Petersson metric on $V_{1}^{k}$ is defined
as 
$\omega_{WP}=-\sqrt{-1}\partial\bar\partial\log(\Omega,\bar\Omega)$ where 
$(\Omega,\bar\Omega)>0$ on $V_{1}^{k}$. We also assume that 
$\Omega=e^{\frac{\sqrt{-1}}{2\pi}\sum_{1}^{k}N_{i}\log\frac{1}{z_i}}A(z)$ 
and 
$\Omega_0=e^{\frac{\sqrt{-1}}{2\pi}\sum_{1}^{k}N_{i}\log\frac{1}{z_i}}A_0$, 
where $A_0=A(0)$. 

The following lemma is one of the key parts in the proof of
Theorem~\ref{thm11}.
\begin{lemma}\label{reduce}
There exist a universal constant $\delta >0$ which only depends on
$\Omega$ such that
$\log(\Omega,\bar\Omega)$  is integrable on $V_{\delta}^{k}$ with respect
to the standard Poincar\'e  metric $\omega_{P}$ on $V_{1}^{k}$.
\end{lemma}
{\bf Proof.}
The proof depends on the convexity of $\log(\Omega,\bar\Omega)$. 
Without loss of generality we can assume the convergence radius of $A(z)$
is $1$. Obviously, we have the upper bound\begin{eqnarray}\label{m10}
\log(\Omega,\bar\Omega) \leq c+\log\Pi_{j=1}^{k}(\log\frac{1}{r_j})^{n}
\end{eqnarray}
where $r_j$ is the Euclidean norm of $z_j$. Clearly
the right hand side of the above formula is integrable 
with respect to the Poincar\'e metric $\omega_{P}$ which
implies that we only need to check that the integration of
$\log(\Omega,\bar\Omega)$ has a lower bound. 
Let $\theta_j$ be the argument of $z_j$ and set 
\begin{eqnarray}\label{m20}
p(r_1,\cdots,r_m)=\int_{0}^{2\pi}\cdots\int_{0}^{2\pi}
\log(\Omega,\bar\Omega)\ d\theta_1\cdots d\theta_m.
\end{eqnarray} 
Since $-\sqrt{-1}\partial\bar\partial\log(\Omega,\bar\Omega)>0$ we have 
$\frac{\partial^2 p}{\partial r_{j}^{2}}+\frac{1}{r_j}\frac{\partial p}
{\partial r_j}<0$ for each $1\leq j \leq m$ which is equivalent to 
\begin{eqnarray}\label{m30}
\frac{\partial (r_j\frac{\partial p}{\partial r_j})}{\partial r_j}<0.
\end{eqnarray}
In order to prove $\log(\Omega,\bar\Omega) \in
L^{1}(V_{\delta}^{k},\omega_{P})$ we just need to check that 
\begin{eqnarray}\label{m40}
\int_{0}^{\delta}\cdots\int_{0}^{\delta}p(r_1,\cdots,r_m)
\frac{r_{k+1}\cdots r_m}{r_{1}
(\log\frac{1}{r_1})^2\cdots r_{k}(\log\frac{1}{r_k})^2}
\ dr_1\cdots dr_m > -\infty.
\end{eqnarray}

We prove this using mathematical induction on
$m$. If the dimension $m$ of the moduli 
space is $1$,
then 
$k$ must be $1$. 
By a theorem of Wang~\cite{wang-1}, we know that the
leading term in $(\Omega,\bar\Omega)$ is $c(\log\frac{1}{r_1})^l$ with $l
\geq 1$ if the
Weil-Petersson metric is complete at $0$ where $c$ is a positive constant. 
So there is a constant $\delta >0$ such that when $r<\delta$,
$(\Omega,\bar\Omega) \geq \frac{c}{2}(\log\frac{1}{r_1})^l$ which implies
that $p(r_1) \geq 2\pi \log\frac{c}{2}+2\pi
l\log\log\frac{1}{r_1}$. So 
\begin{eqnarray*}
\int_{0}^{\delta}p(r_1)\frac{1}{r_1(\log\frac{1}{r_1})^2}\ dr_1
&\geq& 2\pi\log\frac{c}{2}\int_{0}^{\delta}\frac{1}
{r_1(\log\frac{1}{r_1})^2}\ dr_1 
+2\pi l \int_{0}^{\delta}\frac{\log\log\frac{1}{r_1}}
{r_1(\log\frac{1}{r_1})^2}\ dr_1\\
&=&-\frac{2\pi\log\frac{c}{2}}{\log\delta}-\frac{2\pi l}{\log\delta}
(\log\log\frac{1}{\delta}+1)>-\infty.
\end{eqnarray*}

If the Weil-Petersson metric is incomplete at $0$, the leading term of 
$(\Omega,\bar\Omega)$ is a positive constant $c$.
So we can find $\delta>0$ such that $p(r_1) \geq c_1$ for some constant
$c_1$ when $r < \delta$. This implies 
\begin{eqnarray*}
\int_{0}^{\delta}p(r_1)\frac{1}{r_1(\log\frac{1}{r_1})^2}\ dr_1
\geq c_1 \int_{0}^{\delta}\frac{1}{r_1(\log\frac{1}{r_1})^2}\ dr_1
=\frac{c_1}{\log\frac{1}{\delta}}>-\infty.
\end{eqnarray*}
Now
we assume that when $m \leq s-1$ and $1\leq k\leq s-1$, the
inequality
\eqref{m40} hold. Let $m=s$. We first fix $r_2,\cdots,r_m$. From
\eqref{m30} we have $\frac{\partial (r_1\frac{\partial p}{\partial
r_1})}{\partial r_1}<0$ which implies $r_1\frac{\partial p}{\partial r_1}$
decreases  as $r_1$ increases. Like the argument above, the leading term
of $p$ in
$r_1$ is either $\log c+l\log\log\frac{1}{r_1}$ or $\log c$ 
where $c$ is a positive function of 
$r_2,\cdots, r_m$ and $l \geq 1$ is a positive integer. In either cases we
have $\lim_{r_1\to 0} r_1\frac{\partial p}{\partial r_1}=0$. So we know
that 
$r_1\frac{\partial p}{\partial r_1}<0$ when $r_1>0$ which implies 
$p$ is decreasing as $r_1$ is increasing. So for any $\delta_1>0$ small
enough, we have 
\begin{eqnarray}\label{m50}
\int_{0}^{\delta_1} 
\frac{p(r_1,\cdots,r_m )}{r_1(\log\frac{1}{r_1})^2}\ dr_1 \geq 
\int_{0}^{\delta_1} \frac{p(\delta_1,r_2,\cdots,r_m )}
{r_1(\log\frac{1}{r_1})^2}\
dr_1=-\frac{1}{\log\delta_1}p(\delta_1,r_2,\cdots,r_m ).
\end{eqnarray}
Now we go back to $(\Omega,\bar\Omega)$. We fix a $z_1=w \in
\Delta^{\ast}$. Let $\widetilde{\Omega}=\Omega(w,z_2,\cdots,z_m)$. Then
$-\sqrt{-1}\tilde\partial\bar{\tilde\partial}\log(\widetilde{\Omega},
\bar{\widetilde{\Omega}})$
gives the Weil-Petersson metric on the slice $\tilde U = \{ w \} \times
(\Delta^{\ast})^{k-1}\times\Delta^{m-k}$. If $k>1$, 
by the induction assumption, we
know that $\log(\widetilde{\Omega},\bar{\widetilde{\Omega}})$ is
integrable on $\tilde U$ with respect to the Poincar\'e metric for some 
$\delta>0$ and the integration depends on $w$ continuously. This clearly
implies 
\begin{eqnarray}\label{m54}
\int_{0}^{2\pi}
\int_{(\Delta_{\delta}^{\ast})^{k-1}\times\Delta_{\delta}^{m-k}}
\log(\widetilde{\Omega},\bar{\widetilde{\Omega}}) \omega_{P}^{m-1}\
d\theta_1>-\infty.
\end{eqnarray}
Thus we have 
\begin{eqnarray}\label{m56}
\int_{0}^{\delta}\cdots\int_{0}^{\delta}
-\frac{1}{\log\delta_1}p(\delta_1,r_2,\cdots,r_m )
\frac{r_{k+1}\cdots r_m}{r_{1}
(\log\frac{1}{r_1})^2\cdots r_{k}(\log\frac{1}{r_k})^2}
\ dr_2\cdots dr_m > -\infty.
\end{eqnarray}
So $\log(\Omega,\bar\Omega)$ is integrable over $V_{\delta}^{k}$ for some
$\delta>0$ with respect to the Poincar\'e metric.
If $k=1$, then $\tilde U=\{ w \} \times \Delta^{m-1}$ and $\log (
\widetilde\Omega,\bar{\widetilde\Omega})$ is smooth on $\tilde U$. Clearly
\eqref{m54}, \eqref{m56} still hold for some $\delta>0$. So 
$\log(\Omega,\bar\Omega)$ is integrable with respect to the
 Poincar\'e
metric. 

\qed

{\bf Proof of Theorem~\ref{thm11}.}
By the Nilpotent Orbit theorem we know that the
Hodge bundle
$\underline{F}^n$ over ${M}$ can be extended to
$\bar{{M}}$ smoothly. We denoted the extended bundle by
$\underline{E}^n$. Now we put a Hermitian metric $g$ on
$\underline{E}^n$. Take an open cover of
$\bar{{M}}$ like we did in \S \ref{cut}. On each $U_{\alpha}$,
let $g_{\alpha}$ be a representation of  $g$. The local potential of the
Weil-Petersson metric on
$U_{\alpha}$ is given by $h_{\alpha}=(\Omega,\bar\Omega)$. Let  
$f=\frac{h_{\alpha}}{g_{\alpha}}$. It is clear that $f$ is a global
function on ${M}$ although $g_{\alpha}$ and $h_{\alpha}$ are only
locally defined. Let $\tilde\omega$ be the  curvature form of
the metric $g$. Then 
$\omega_{WP}=-\bb\partial\bar\partial\log f+\tilde\omega$.
Since we know that
$\tilde\omega$ is the first Chern class of the line bundle
$\underline{E}^n$ over the K\"ahler manifold $\bar{{M}}$, we know
that
\begin{eqnarray*}
\int_{{M}}\tilde\omega^{m}
=\int_{\bar{{M}}}\tilde\omega^{m}
\end{eqnarray*} 
is an integer. We  need to prove that 
$\int_{{M}}\tilde\omega^{m}
=\int_{{M}}\omega_{WP}^{m}$. Using Schwarz-Yau lemma (Proposition~\ref{SY})
 like we
did in \cite{LS-1} we know that $\omega_{WP} \leq c\omega_{P}$ where 
$\omega_{P}$ is the asymptotic Poincar\'e metric we constructed
in section
\ref{cut}. Since $\tilde\omega$ is the Ricci form of a line bundle over  the compact
manifold $\bar{{M}}$, it is bounded. Thus we can find a constant $c$ such that 
$-c\omega_{P}\leq \tilde\omega \leq c\omega_{P}$. 

We check the integrability of $\log f$  on ${M}$ with respect
to the asymptotic metric $\omega_{P}$. Let $\psi_{\alpha}$ be
a partition of unity subordinated to the cover
$\{U_{\alpha}\}$. For a chart $U_{\alpha}$ of $\bar{{M}}$, if 
$U_{\alpha} \cap \bar{{M}}\setminus{M}=\emptyset$ then 
$\psi_{\alpha}\log f$ is bounded on $U_{\alpha}$. So 
$\int_{U_{\alpha}}\psi_{\alpha}\log f\ \omega_{P}^{m}$ is
finite. If 
$U_{\alpha} \cap \bar{{M}}\setminus{M}\ne\emptyset$, we know
that $\psi_{\alpha}\log f=\psi_{\alpha}\log(\Omega,\bar\Omega)-
\psi_{\alpha}\log g_{\alpha}$. Clearly $\psi_{\alpha}\log g_{\alpha}$ is
bounded on $U_{\alpha}$. By Lemma \ref{reduce} we know that 
$\int_{U_{\alpha}}\psi_{\alpha}\log(\Omega,\bar\Omega)
\ \omega_{P}^{m}$ is also finite. This implies $\log f$ is
integrable. 

Pick an $\epsilon>0$ small. Let $\rho_{\epsilon}$ be the cut-off function
we constructed in section \ref{cut}. We have 
\begin{eqnarray}\label{m60}
\int_{{M}}(\rho_{\epsilon}\tilde\omega^{m}
-\rho_{\epsilon}\omega_{WP}^{m})=\sum_{j=0}^{m-1}\int_{{M}}
\rho_{\epsilon}\tilde\omega^{j}\wedge\omega_{WP}^{m-j-1}\wedge
(\tilde\omega-\omega_{WP}). 
\end{eqnarray}
For each $0\leq j \leq m-1$, we have 
\begin{align}\label{m70}
\begin{split}
& \int_{{M}}
\rho_{\epsilon}\tilde\omega^{j}\wedge\omega_{WP}^{m-j-1}\wedge
(\tilde\omega-\omega_{WP})= \int_{{M}}
\rho_{\epsilon}\tilde\omega^{j}\wedge\omega_{WP}^{m-j-1}\wedge 
\partial\bar\partial\log f \\
= &  \int_{{M}}\log f \partial\bar\partial\rho_{\epsilon}
\wedge \tilde\omega^{j}\wedge\omega_{WP}^{m-j-1}
=\int_{\rm supp(1-\rho_{\epsilon})}\log f
\partial\bar\partial\rho_{\epsilon}
\wedge \tilde\omega^{j}\wedge\omega_{WP}^{m-j-1}\\
=& \sum_{\alpha}
\int_{\rm supp(1-\rho_{\epsilon}) \cap U_{\alpha}}\psi_{\alpha}\log f
\partial\bar\partial\rho_{\epsilon}
\wedge \tilde\omega^{j}\wedge\omega_{WP}^{m-j-1}
\end{split}
\end{align}
where the sum over $\alpha$ is a finite sum. 
Also, on each $U_{\alpha}$, $\log f$ is bounded above by a
positive function 
$c+\sum_{1}^{l_{\alpha}}\log\log\frac{1}{r_j}$ which   is
integrable with respect to the local Poincar\'e metric. 
By Theorem
\ref{thm41} and the fact that 
$\tilde \omega+\omega_{WP}\leq c\omega_P$ (which follows from
Proposition~\ref{SY}),
 we know that 
\begin{align}\label{m80} 
\begin{split}
&\left |
\int_{\rm supp(1-\rho_{\epsilon}) \cap U_{\alpha}}\psi_{\alpha}\log f
\partial\bar\partial\rho_{\epsilon}
\wedge \tilde\omega^{j}\wedge\omega_{WP}^{m-j-1}\right |\\
\leq & \left |
\int_{\rm supp(1-\rho_{\epsilon}) \cap U_{\alpha}}\psi_{\alpha}
(\log f+c+\sum_{1}^{l_{\alpha}}\log\log\frac{1}{r_j})
\partial\bar\partial\rho_{\epsilon}
\wedge \tilde\omega^{j}\wedge\omega_{WP}^{m-j-1}\right |\\
&+\left |
\int_{\rm supp(1-\rho_{\epsilon}) \cap U_{\alpha}}\psi_{\alpha}
(c+\sum_{1}^{l_{\alpha}}\log\log\frac{1}{r_j})
\partial\bar\partial\rho_{\epsilon}
\wedge \tilde\omega^{j}\wedge\omega_{WP}^{m-j-1}\right |\\
\leq & c_1\left |
\int_{\rm supp(1-\rho_{\epsilon}) \cap U_{\alpha}}\psi_{\alpha}
(\log f+c+\sum_{1}^{l_{\alpha}}\log\log\frac{1}{r_j})
\omega_{P}^{m}\right |\\
&+c_1\left |
\int_{\rm supp(1-\rho_{\epsilon}) \cap U_{\alpha}}\psi_{\alpha}
(c+\sum_{1}^{l_{\alpha}}\log\log\frac{1}{r_j})
\omega_{P}^{m}\right|.
\end{split}
\end{align}
The above expression converges to $0$ as $\epsilon \to 0$,
because
$\omega_{P}$ has finite volume and the measure of ${\rm
supp}(1-\rho_{\epsilon})$ goes to
$0$ as $\epsilon \to 0$ and both 
$\log f+c+\sum_{1}^{l_{\alpha}}\log\log\frac{1}{r_j}$ and 
$c+\sum_{1}^{l_{\alpha}}\log\log\frac{1}{r_j}$ 
are non-negative and  integrable with respect to the
metric $\omega_{P}$. Combining \eqref{m60}, \eqref{m70} and \eqref{m80}
we have 
\begin{eqnarray*}
\lim_{\epsilon\to 0}\int_{{M}}(\rho_{\epsilon}\tilde\omega^{m}
-\rho_{\epsilon}\omega_{WP}^{m})=0.
\end{eqnarray*}
Thus
\begin{eqnarray*}
\int_{{M}}(\tilde\omega^{m}-\omega_{WP}^{m})=\lim_{\epsilon\to 0}
\int_{{M}}(\rho_{\epsilon}\tilde\omega^{m}
-\rho_{\epsilon}\omega_{WP}^{m})=0.
\end{eqnarray*} 
This finishes the
proof.  In general,
if the Nilpotent operators $\{N_i\}$ are not unipotent, then 
by Lemma~\ref{lem41}, we know that the volume is at least 
a rational number.

\qed

\section{First Chern class}\label{chernclass}

Let $ M$ be a Weil-Petersson variety of
dimension $m$ and let $\omega_{WP}$ be the K\"ahler
form of the Weil-Petersson metric. 
Like in the previous sections, we let $\bar M$ be 
the compactification of $M$ such that $Y=\bar
M\backslash M$ is a divisor of normal crossings.

The main result of
this section is the following
\begin{theorem}\label{c1}
Let $\bar X \subset \bar{M}$ be a subvariety of
dimension
$q$. Let $X=\bar X\cap M$. Then
\begin{eqnarray*}
\int_{X}(Ric(\omega_{WP}))^{k}\wedge \omega_{WP}^l \in
\mathbb{Q}
\end{eqnarray*}
for $k+l=q$, where the integration is on the
smooth part of the variety $X$.
\end{theorem}
{\bf Proof.} Without loss generality,
we assume that $q=m$ and $X=M$.
Let $L_0=K_M^{-1}$ be the anti-canonical
line bundle of  $M$. By Lemma~\ref{lem41},
we assume that all elements of the
 monodromy group are unipotent. 
By the remark in \S \ref{ext} we know
that the tangent bundle of $ M$, as a quotient
of the Hodge bundles, can be extended to the
compactification $\bar{ M}$. This implies the
$L_0$ can be extended to $\bar{L_0}$ over $\bar{
M}$.

Now let $Y=\bar{ M}\setminus {M}$ be the  divisor of
$\bar{ M}$ of normal crossings. Let $L_j$ be the
line bundle corresponding to $Y_j$ for $j=1,\cdots,p$. We write
down the transition functions of $L_j$ explicitly. Recall that in
Section \ref{ext}, we constructed an open cover $\{ U_{\alpha,0}
\} \cup \{ U_{\alpha,i} \} \cup \{ V_{\beta} \}$ of $\bar{
M}$. Denote this cover by $\mathfrak U$. We knew that there is
a neighborhood $S$ of $Y$ in $\bar{ M}$ such that $\bar{S}
\subset \bigcup U_{\alpha,0}$. Let
$\widetilde{U}_{\alpha,j}=U_{\alpha,j}\setminus \bar{S}$. Then $\{
U_{\alpha,0} \} \cup \{ \widetilde{U}_{\alpha,i} \} \cup \{
V_{\beta} \}$ is also an open cover of $\bar{ M}$, denoted
by $\widetilde{\mathfrak  U}$. For each $j=1,\cdots,p$, let
$I_j$ be the index set such that $\alpha \in I_j$ if and only if
$U_{\alpha,0}\cap Y_j\ne\emptyset$. From \S \ref{ext} we can
assume that for each $\alpha \in I_j$, $U_{\alpha,0}$ has
coordinates $z^\alpha_1,\cdots,z^\alpha_m$ and $Y_j \cap
U_{\alpha,0}$ is given by $z^\alpha_1=0$.

Now we define the transition functions. By reordering the elements
in $\widetilde{\mathfrak U}$, we assume that
$\tilde{\frk{{U}}}= \{ U_\alpha \}_{\alpha \in I_j} \cup \{
W_\gamma \}$,
where $U_\alpha\cap Y_j\neq\emptyset$ and $W_\gamma 
\cap Y_j=\emptyset$. When $\alpha,\beta \in I_j$ and
$U_{\alpha,0}\cap U_{\beta,0}\ne\emptyset$, the transition
functions are defined to be
$g_{_{U_{\alpha,0},U_{\beta,0}}}=\frac{z^\alpha_1}{z^\beta_1}$.
For $\alpha \in I_j$, define
$g_{_{U_{\alpha,0},W_\gamma}}=z^\alpha_1$ for each $W_\gamma \in
\tilde{
\frk{ U}}$ with $U_{\alpha,0}\cap
W_\gamma\ne\emptyset$. Finally,
for $W_\gamma, W_\delta \in \tilde{\frk{{U}}}$ with
$W_\gamma \cap W_\delta\ne \emptyset$, define
$g_{_{W_\gamma,W_\delta}}=1$. One can easily check that these
transition functions define the line bundle $L_j$.

Now we define the degeneration order of the volume form of the
Weil-Petersson metric along each hypersurface $Y_j$. We will need
the following lemma which is proven in \cite{LS-1}. For
completeness, we include the proof here. In the following, we use
$\partial=\partial_{z}$ and $\bar\partial=\bar{\partial_{z}}$.
\begin{lemma}\label{pos}
Let $f:\mathbb{C} \to R$ be a degree $k$ homogeneous polynomial of
$z$ and $\bar z$. Assume that $f(z,\bar z) \geq 0$ and is not
identically $0$. If there is a constant $c_1>0$ such that for
every $z$ with $f(z,\bar z)\ne 0$, we have
\[
-\frac{c_1}{r^2(\log\frac{1}{r})^2} \leq -
\partial \bar\partial 
\log f \leq \frac{c_1}{r^2(\log\frac{1}{r})^2}
\]
where $r=|z|$, then $k=2l$ is an even integer and $f(z,\bar
z)=cr^k$ where $c$ is a positive constant.
\end{lemma}
{\bf Proof.} Let $\theta$ be the argument of $z$. Since $f$ is a
homogeneous polynomial which is not identically $0$, we can easily
see that the set $\{ z \in S^1 \mid f(z,\bar z)=0 \}$ is a
$0$-dimensional real analytic variety and thus is a discrete set
of $S^1$. So there are only finite many rays from the origin where
$f$ vanishes.

If $k$ is odd, since each term of $f$ has the form
$pz^l\bar{z}^{k-l}$, we know that
\[
\int_{0}^{2\pi} f(e^{i\theta},e^{-i\theta}) \ d\theta=0,
\]
which contradicts the fact that $f\geq 0$ and
$f(e^{i\theta},e^{-i\theta})$ vanish only for finitely many
$\theta$. So $k$ has to be an even integer.

We have
\[
-\partial\bar\partial\log f=\frac{\partial f\bar\partial f
-f\partial\bar\partial f}{f^2}.
\]
Clearly if $\partial f\bar\partial f-f\partial\bar\partial f$ is
not identically $0$, then it is a degree $2k-2$ homogeneous
polynomial. Consequently, $-\partial\bar\partial\log f$ is of
order $\frac{1}{r^2}$ which contradicts to the assumption. Thus we
have
\begin{eqnarray}\label{iden0}
\partial f\bar\partial
f-f\partial\bar\partial f=0.
\end{eqnarray}

It is also clear that for a homogeneous polynomial of $z$ and
$\bar z$, if it is identically $0$, then all of its coefficients
are $0$.  Now we use induction on $l$. When $l=1$, we have
$f=a_0z^2+a_1z\bar z+a_2\bar{z}^2$. From \eqref{iden0} we know
that $a_0a_1=a_0a_2=a_1a_2=0$. If either $a_0\ne 0$ or $a_2\ne 0$
we know that the rest of the coefficients are $0$ which implies
that $f$ is not real. So $a_0=a_2=0$. Since $f$ is real and
non-negative, we know that $a_1 >0$. The lemma holds. Assuming
that the lemma hold when $l \leq p-1$, 
we consider the case $l=p$.
Assume that $f=a_0z^k+a_1z^{k-1}\bar{z}+\cdots+a_k\bar{z}^k$ where
$k=2p$. If $a_0\ne 0$, then we assume $i=\underset{j>0,a_j\ne
0}{\min} j$. Consider the term $z^{2k-i-1}\bar{z}^{i-1}$ in
$\partial f\bar\partial f-f\partial\bar\partial f$. We have
$kia_0a_i-i(k-i)a_0a_i=0$ which implies $a_i=0$ since $a_0\ne 0$.
This means $f=a_0z^k$ which contradicts the fact that $f$ is real.
So $a_0=0$. Similarly we can prove that $a_k=0$. So
$f=z\bar{z}f_0$ where $f_0$ is a homogenous polynomial of $z,\bar
z$ of degree $k-2$. $f_0$ satisfies the conditions of the lemma.
By the assumption, $f_0=cz^{p-1}\bar{z}^{p-1}$. So
$f=cz^{p}\bar{z}^{p}=cr^k$.

\qed

Pick a point $q \in Y_j$ and assume $q \in U=U_{\alpha,0}$ with
the local coordinates
$z=(z_1,\cdots,z_m)=z_\alpha=(z^\alpha_1,\cdots,z^\alpha_m)$ on
$U$. Let $g_{i\bar j}$ be the Weil-Petersson metric. In the
following, we will use $r_j$ and $\theta_j$ to denote the
Euclidean norm and argument of $z^j$. We assume $q \in Y_j
\setminus \bigcup_{k\ne j}Y_k$. Roughly speaking, the
degenerate order of $\omega_{WP}^m$ along $Y_j$ measures 
the rate of the blow-up or degeneration of the volume form
of the Weil-Petersson metric.
We need to analyze the asymptotic behavior of $\omega_{WP}^m$
when $z \to q\in Y_j$. 
\begin{lemma}
Assume $Y_j\cap U$ is given by $z_1=0$ and $q
\in Y_j \setminus \bigcup_{k\ne j}Y_k$. Then when we expand
$(\Omega,\bar\Omega)^{2m}\det(g_{i\bar j})$, the leading term
in $z_1$ has form
$A_k(z',\bar{z'})r_{1}^{k-2} (\log\frac{1}{r_1^2})^l$ where
$z'=(z_2,\cdots,z_m)$.
\end{lemma}
{\bf Proof.} Near $q$, we know that the local holomorphic section
$\Omega$ of the first Hodge bundle $F^n$ can be written as
$\Omega=e^{\frac{\sqrt{-1}}{2\pi}N_1\log\frac{1}{z_1}}
A(z_1,\cdots,z_m)$ where $A$ is a holomorphic function. Let
$z'=(z_2,\cdots,z_m)$. We can expand $A$ as a power series of
$z_1,\cdots,z_m$ and assume that the convergent radius is $1$. Let
$g=\det(g_{i\bar j})$. For each $i=1,\cdots,m$, let
$\Omega_i=e^{\frac{\sqrt{-1}}{2\pi}N_1\log\frac{1}{z_1}}
\partial_i A$ and let
$\widetilde{\Omega}_1=z_1\Omega_1-\bb N_1\Omega$. Let
\begin{enumerate}
\item $\widetilde{g_{1\bar 1}}=(\widetilde{\Omega}_1,\bar\Omega)
(\Omega,\bar{\widetilde{\Omega}_1})-(\Omega,\bar\Omega)
(\widetilde{\Omega}_1,\bar{\widetilde{\Omega}_1})$;
\item $\widetilde{g_{1\bar j}}=(\widetilde{\Omega}_1,\bar\Omega)
(\Omega,\bar{\Omega_j})-(\Omega,\bar\Omega)
(\widetilde{\Omega}_1,\bar{\Omega_j})$ for $2 \leq j \leq m$;
\item $\widetilde{g_{i\bar 1}}=(\Omega_i,\bar\Omega)
(\Omega,\bar{\widetilde{\Omega}_1})-(\Omega,\bar\Omega)
(\Omega_i,\bar{\widetilde{\Omega}_1})$ for $2 \leq i \leq m$;
\item $\widetilde{g_{i\bar j}}=(\Omega_i,\bar\Omega)
(\Omega,\bar{\Omega_j})-(\Omega,\bar\Omega)
(\Omega_i,\bar{\Omega_j})$ for $2\leq i,j \leq m$.
\end{enumerate}
Let $\widetilde{g}=\det(\widetilde{g_{i\bar j}})$. We have
$g=\frac{1}{(\Omega,\bar\Omega)^{2m}}
\frac{1}{r_1^2}\widetilde{g}$. Since $g>0$ we know that
$\widetilde{g}>0$. We expand $A$ with respect to $z_1$ so that the
coefficient of each term $z_1^s$ is an analytic function of
$z'$ and $\bar{z'}$. By monodromy theorem,
$N_1^{n+1}=0$, where $N_1$ is the nilpotent operator. Also,
since
$N_1$ is nilpotent,
we know that each term in
$\widetilde{g}$ has the form
$A_{s,t,l}z_1^s\bar{z_1}^t(\log\frac{1}{r_1^2})^l$ where
$A_{s,t,l}=A_{s,t,l}(z',\bar{z}')$ and $l \leq mn$. Define the
degree of such a term by $s+t-\frac{l}{mn+1}$. Assume the lowest
degree of terms in the expansion of $\widetilde{g}$ is
$k-\frac{l}{mn+1}$. Collect all the terms of degree
$k-\frac{l}{mn+1}$ having the form $\sum_{s=0}^k
A_{s,k-s,l}z_1^s\bar{z_1}^{k-s}(\log\frac{1}{r_1^2})^l$. Let
$f(z_1,\bar{z_1})=\sum_{s=0}^k A_{s,k-s,l}z_1^s\bar{z_1}^{k-s}$.
We know that, except for a set of lower dimension,
$f(z_1,\bar{z_1}) \ne 0$ and $f(z_1,\bar{z_1})$ is a homogeneous
polynomial of $z_1$ and $\bar{z_1}$ whose coefficients are functions
of $z'$ and $\bar{z'}$ . For each fixed $z'$ with
$f(z_1,\bar{z_1}) \ne 0$, since
$f(z_1,\bar{z_1})(\log\frac{1}{r_1^2})^l$ is the leading term in
the expansion of $\widetilde{g}$, we know that $f(z_1,\bar{z_1})
\geq 0$ because $\widetilde{g} >0$. We call a point $z'$ a generic
point in the first direction if for this $z'$, we have
$f(z_1,\bar{z_1})
>0$. For a generic point $z'$, we know that
$\widetilde{g}=f(z_1,\bar{z_1}) (\log\frac{1}{r_1^2})^l+g_0$ where
each term in $g_0$ has degree higher than $k-\frac{l}{mn+1}$. 
Finally we
have
\begin{align*}
&\qquad -\bb\partial\bar\partial\log\widetilde{g}=
-\bb\partial\bar\partial\log
(r_1^2(\Omega,\bar\Omega)^{2m}g)\\&
=-2m\bb\partial\bar\partial\log(\Omega,\bar\Omega)-\bb
\partial\bar\partial\log g=2m\omega_{WP}+Ric(\omega_{WP}).
\end{align*}
By the Strominger's formula we know that there is a positive
constant $c$ such that
\[
-c\omega_{P} \leq Ric(\omega_{WP}) \leq c\omega_P,
\]
where $\omega_P$ is the K\"ahler form of the Poincar\'e metric.
Finally by using Hodge metric and  Schwarz-Yau lemma (Proposition~\ref{SY}) we proved in
\cite{LS-1} that there is a constant $c$ such that
\[
0<\omega_{WP} \leq c\omega_P.
\]
Combine all these formulae we know that there is a constant $c$
such that
\[
-c\omega_{P} \leq -\bb\partial\bar\partial\log\widetilde{g}
\leq c\omega_P,
\]
which implies
\[
-\frac{c}{r_1^2(\log\frac{1}{r_1})^2} \leq
-\partial_1\bar{\partial_1}
\log \widetilde{g} \leq \frac{c}{r_1^2(\log\frac{1}{r_1})^2}.
\]
Since
\[
-\partial_1\bar{\partial_1}
\log
\widetilde{g}=\frac{\partial_1\widetilde{g}
\bar{\partial_1}\widetilde{g}-\widetilde{g}
\partial_1\bar{\partial_1}\widetilde{g}
}{\widetilde{g}^2}
\]
and the leading term in the numerator is
$(\partial_1 f\bar{\partial_1} f-f\partial_1\bar{\partial_1}f)
(\log\frac{1}{r_1^2})^{2l}$, the leading term in the denominator
is $f^2(\log\frac{1}{r_1^2})^{2l}$, we know that, for generic
$z'$ we have
\[
-\frac{c}{r_1^2(\log\frac{1}{r_1})^2} \leq
-\partial_1\bar{\partial_1}
\log f \leq \frac{c}{r_1^2(\log\frac{1}{r_1})^2}.
\]
Using
the Lemma~\ref{pos}, we have that for generic point $z'$ in the
first direction, $f(z_1,\bar{z_1})=A_k(z',\bar{z'})r_{1}^k$.
This implies that, except a lower dimensional set of $z'$, the
leading term in the expansion of $(\Omega,\bar\Omega)^{2m}
\det(g_{i\bar j})$ is $A_k(z',\bar{z'})r_{1}^{k-2}
(\log\frac{1}{r_1^2})^l$.

\qed

\begin{definition}\label{def666}
The degeneration order of $\omega_{WP}^m$ along
$Y_j$ is $\tau_j=(k-2)/2=\mu_j$, which is an  integer
by Lemma~\ref{pos}. 
\end{definition}

\begin{lemma}
The degeneration order is well-defined and is constant except
for a lower dimensional set.
\end{lemma}
{\bf Proof.}
Assume there is another chart
$U_\beta$ with $U_\alpha \cap U_\beta \cap Y_j \ne \emptyset$.
We can also assume that $U_\beta \cap Y_j$ is given by
$z_1^\beta=0$. From the definition we can see that the
degeneration order $\tau_j$ is the smallest number
$\tau$ such that, for a generic point $z'$ in the first
direction, $\liminf_{z_1 \to 0}r_1^\tau(\Omega,\bar\Omega)^{2m}
\det(g_{i\bar j})>0$. Since
\[
\omega_{WP}^m=(\frac{\sqrt{-1}}{2\pi})^m
m!\det((g_\alpha)_{i\bar j}) dz_1^\alpha\wedge
d\bar{z_1^\alpha}\cdots  dz_m^\alpha\wedge d\bar{z_m^\alpha}
\]
is a global form and we can choose the same $\Omega$ on
$U_\alpha\cap U_\beta$, we have
\[
(r_1^\alpha)^\tau \det((g_\alpha)_{i\bar j})
=\bigg(\frac{r_1^\alpha}{r_1^\beta}\bigg)^\tau
\left | \det(\frac{\partial z_i^\beta}{\partial
z_j^\alpha})\right | ^2 (r_1^\beta)^\tau\det((g_\beta)_{i\bar
j}).
\]
However, we know that $\bigg(\frac{r_1^\alpha}{r_1^\beta}\bigg)^\tau
\left | \det(\frac{\partial z_i^\alpha}{\partial
z_j^\beta})\right | ^2$ is bounded and bounded below from $0$.
When we choose $U_\beta=U_\alpha$ but with a different
coordinate system, the above argument implies the degeneration
order is independent of the choice of local coordinates. For
general $U_\alpha$ and $U_\beta$, this implies that the
order is constant along $Y_j$ except for a lower dimensional
set.

\qed

Now we prove that $\int_{{M}}
(Ric(\omega_{WP}))^{s}\wedge \omega_{WP}^{m-s} \in \mathbb Q$.
We first extend the first Hodge bundle $F^n$ to the
compactification $\bar{ M}$ in the way described in
\S \ref{ext}. We put a smooth Hermitian metric $h$ on this
extended bundle and denote its K\"ahler form by $\omega$. Recall
that
we use $L_0$ to denote the anti-canonical bundle
of $ M$ and use $\bar{L_0}$ to denote the preferred 
extension of $L_0$ to
$\bar{ M}$. We put a smooth Hermitian metric $h_0$ on
$\bar{L_0}$. Finally for each line bundle
$L_j$ corresponding to the hypersurface $Y_j$, we put a smooth
Hermitian metric $h_j$ on it and denote its curvature form by
$\omega_j$.

We first check that, for all $1\leq s \leq m$, $\int_{
M}(Ric(\omega_{WP}))^s\wedge \omega_{WP}^{m-s} \in \mathbb{Q}$ is
equivalent to
\begin{eqnarray}\label{s-10}
\int_{ M}(Ric(\omega_{WP})+2m\omega_{WP})^s\wedge
\omega_{WP}^{m-s} \in \mathbb{Q}
\end{eqnarray}
for all $1\leq s \leq m$. Clearly we know
$(Ric(\omega_{WP})+2m\omega_{WP})^s\wedge
\omega_{WP}^{m-s}$ is an integer linear combination of terms
like $(Ric(\omega_{WP}))^i\wedge \omega_{WP}^{m-i}$ and
$ \omega_{WP}^{m}$. Also, $(Ric(\omega_{WP}))^s\wedge
\omega_{WP}^{m-s}$ is an integer linear combination of terms
like $(Ric(\omega_{WP})+2m\omega_{WP})^j\wedge
\omega_{WP}^{m-j}$ and $\omega_{WP}^{m}$. Since the
Weil-Petersson volume is an integer, it is clear that the above
argument is true.

The reason we use the expression in~\eqref{s-10} is that
${\rm Ric} (\omega_{WP})+2m\omega_{WP}$ is a nonnegative
form.

Realizing that the Ricci curvature of the Weil-Petersson metric
is bounded above and below by a constant multiple of the
asymptotic Poincar\'e metric, we can use the same proof in
Section ~\ref{vol} to derive that
\begin{eqnarray}\label{s0}
\int_{ M}(Ric(\omega_{WP})+2m\omega_{WP})^s\wedge
\omega_{WP}^{m-s}=
\int_{ M}(Ric(\omega_{WP})+2m\omega_{WP})^s\wedge
\omega^{m-s}_0,
\end{eqnarray}
where $\omega_0$ is the curvature form of the line bundle
$\bar L_0$.

Recall that we use $\mu_j$ to denote the degeneration order of
$\omega_{WP}^m$ along the hypersurface $Y_j$. Let
$\widetilde{\omega}=-\partial\bar\partial\log h_0
+\sum_{j=1}^p\mu_j\omega_j+2m\omega_{WP}$. Then
\begin{eqnarray}\label{s10}
\widetilde{\omega}=c_1
(\bar{L_0}\oplus L_1^{\mu_1}\oplus\cdots\oplus
L_p^{\mu_p}\oplus(\bar{F^n})^{2m}).
\end{eqnarray}

For each line bundle $L_j$, on a chart $U_\alpha$
at
$Y_j$, we assume the metric $h_j$ is given by $h_j^\alpha$. From
the transition functions that define $L_j$ we can see that
$h_j^\alpha|z_1^\alpha|^2$ is a global function on
$\bar{ M}$. Denote  $f_j=(h_j^\alpha |z_1^\alpha|^2)^{\mu_j}$.
By reordering coordinates, we know that
\begin{eqnarray}\label{s20}
\widetilde{f}=\frac{\omega_{WP}^m}{h_0f_1\cdots f_p}
\end{eqnarray}
is a global function on $ M$. Also, from \S
\ref{vol}, if we use $h^\alpha$ to denote the smooth Hermitian
metric $h$ on $\bar{F^n}$ on a chart $U_\alpha$, we know that
$\bar f=\frac{(\Omega,\bar\Omega)^{2m}}
{(h^\alpha)^{2m}}$ is also a global function.
We have
\begin{align}\label{s30}
\begin{split}
& \int_{ M}
\rho_{\epsilon}(Ric(\omega_{WP})+2m\omega_{WP})^s\wedge
\omega_{0}^{m-s}-\int_{ M}
\rho_{\epsilon}\widetilde{\omega}^s
\wedge \omega_{0}^{m-s}\\
=& \sum_{j=1}^s \int_{ M} \rho_{\epsilon}
(Ric(\omega_{WP})+2m\omega_{WP}-\widetilde{\omega})\wedge
(Ric(\omega_{WP})+2m\omega_{WP})^{j-1}\wedge
\widetilde{\omega}^{s-j}\wedge \omega_{0}^{m-s}\\
=& \sum_{j=1}^s \int_{ M} \rho_{\epsilon}
(-\partial\bar\partial\log\widetilde{f}-\partial\bar\partial
\log f_1)\wedge
(Ric(\omega_{WP})+2m\omega_{WP})^{j-1}\wedge
\widetilde{\omega}^{s-j}\wedge \omega_{0}^{m-s}\\
=& \sum_{j=1}^s \int_{ M} \rho_{\epsilon}
(-\partial\bar\partial\log f)\wedge
(Ric(\omega_{WP})+2m\omega_{WP})^{j-1}\wedge
\widetilde{\omega}^{s-j}\wedge \omega_{0}^{m-s}
\end{split}
\end{align}
where $f=\widetilde{f}f_1$ is a global
positive function on $ M$. Integral by part we have
\begin{align}\label{s40}
\begin{split}
&\int_{ M} \rho_{\epsilon}
(\partial\bar\partial\log f)\wedge
(Ric(\omega_{WP})+2m\omega_{WP})^{j-1}\wedge
\widetilde{\omega}^{s-j}\wedge \omega_{0}^{m-s}\\
=&\int_{ M}\log f
(\partial\bar\partial\rho_{\epsilon})\wedge
(Ric(\omega_{WP})+2m\omega_{WP})^{j-1}\wedge
\widetilde{\omega}^{s-j}\wedge \omega_{0}^{m-s}.
\end{split}
\end{align} 
Like in the case of \S \ref{vol}, we need to prove that $\log
f$ is locally integrable. We use the same notation as we did in
the beginning of \S \ref{vol}.

Consider a chart $U$ at the divisor $Y$. Assume
$U=(\Delta^{\ast})^k\times \Delta^{m-k}$. From our decomposition
$Y=\bigcup Y_j$ we can assume that $Y\cap U=\bigcup_{j=1}^k
(Y_j\cap U)$.
\begin{lemma}
There is a constant $\delta>0$ such that $\log f$ is integrable
on $U_\delta$ with respect to the standard Poincar\'e metric on
$U$.
\end{lemma}
{\bf Proof.}
The proof of this lemma is similar to the proof of Lemma
\ref{reduce} using convexity. On $U$, we have
\begin{align}\label{s50}
\begin{split}
f &=\frac{\omega_{WP}^m}{h_0f_1\cdots f_p}
\frac{(\Omega,\bar\Omega)^{2m}}{h^{2m}} =
\frac{\det{(g_{i\bar j})}(\Omega,\bar\Omega)^{2m}}
{h_0(h_1|r_1|^2)^{\mu_1}\cdots(h_k|r_k|^2)
^{\mu_k}}\frac{1}{f_{k+1}\cdots f_p h^{2m}}\\
&=g_0
\frac{1}{h_0h_1^{\mu_1}\cdots h_k^{\mu_k}
f_{k+1}\cdots f_p h^{2m}}
=g_0g_1,
\end{split}
\end{align}
where $g_0=\frac{\det{(g_{i\bar j})}(\Omega,\bar\Omega)^{2m}}
{r_1^{\tau_1}\cdots r_k^{\tau_k}}$
and
$g_1=\frac{1}{h_0
h_1^{\mu_1}\cdots h_k^{\mu_k}f_{k+1}\cdots f_p h}$. Clearly
$g_1$ is a bounded function and is bounded away from $0$. Since
the Poincar\'e metric on $U$ has finite volume, we know that
$\log g_1$ is locally integrable with respect to the Poincar\'e
metric. Since
\[
-\bb\partial\bar\partial\log f=Ric(\omega_{WP})+
2m \omega_{WP}-\widetilde{\omega} \geq -c\omega_P
\]
and $-\partial\bar\partial\log g_1$ is bounded we have
\begin{eqnarray}\label{s60}
-\bb\partial\bar\partial \log g_0=-\bb\partial\bar\partial\log
f +\bb\partial\bar\partial\log g_1 \geq -c\omega_P.
\end{eqnarray}
From the definition of the degeneration order we know that the
leading term of $g_0$ with respect to $z_j$ for $1 \leq j \leq
k$ is $A_j(z',\bar{z'})(\log\frac{1}{r_j^2})^{l_j}$ with $l_j
\leq mn$ and $z'=(z_1,\cdots,z_{j-1},z_{j+1},\cdots,z_m)$. This
implies that
$g_0
\leq c(\prod_{1}^{k}\log\frac{1}{r_j^2})^{mn}$ which is
integrable with respect to the Poincar\'e metric. So we only
need to prove that
\begin{eqnarray}\label{s70}
\int_U \log g_0 \ \omega_{P}^m > -\infty.
\end{eqnarray}
Set
\begin{eqnarray}\label{s80}
p(r_1,\cdots,r_m)=\int_{0}^{2\pi}\cdots\int_{0}^{2\pi}
\log g_0\ d\theta_1\cdots d\theta_m.
\end{eqnarray}
From \eqref{s60} we can easily see that for $1 \leq j \leq k$,
\begin{eqnarray}\label{s90}
\frac{\partial (r_j\frac{\partial p}{\partial r_j})}{\partial
r_j}\leq \frac{c}{r_j(\log\frac{1}{r_j})^2}.
\end{eqnarray}
Now, for generic $z'$, since the leading term of $g_0$ is
$A_j(z',\bar{z'})(\log\frac{1}{r_j^2})^{l_j}$ where
$A_j(z',\bar{z'})>0$, we know that
\begin{eqnarray}\label{s100}
\lim_{r_j \to 0}r_j\frac{\partial p}{\partial r_j}=0.
\end{eqnarray}
Now integrate both sides of \eqref{s90} from $0$ to $r_j$, using
\eqref{s100}, we have
\begin{eqnarray}\label{s110}
r_j\frac{\partial p}{\partial r_j} \leq
\frac{c}{\log\frac{1}{r_j}}.
\end{eqnarray}
Fix a $\delta>0$ to be chosen, we have that, for generic $z'$
\begin{align*}
\begin{split}
&p(r_1,\cdots,r_{j-1},\delta,r_{j+1},\cdots,r_m)
-p(r_1,\cdots,r_m)\\
\leq & \int_{r_j}^\delta \frac{c}{r_j\log\frac{1}{r_j}}dr_j
=\log\log\frac{1}{r_j}-\log\log\frac{1}{\delta}
\end{split}
\end{align*}
which implies
\begin{eqnarray}\label{s120}
p(r_1,\cdots,r_m) \geq
p(r_1,\cdots,r_{j-1},\delta,r_{j+1},\cdots,r_m)+
\log\log\frac{1}{\delta}-\log\log\frac{1}{r_j}.
\end{eqnarray}
From \eqref{s70} we only need to show that
\begin{align}\label{s130}
\begin{split}
\int_U \log g_0 \ \omega_{P}^m
=\int_{0}^{\delta}\cdots\int_{0}^{\delta}p(r_1,\cdots,r_m)
\frac{r_{k+1}\cdots r_m}{r_{1}
(\log\frac{1}{r_1})^2\cdots r_{k}(\log\frac{1}{r_k})^2}
\ dr_1\cdots dr_m > -\infty.
\end{split}
\end{align} 
We prove \eqref{s130} using induction on $k$. If $k=0$, then $p$
is bounded which implies \eqref{s130} is true.
Assume that for $k \leq l$ the argument is true, consider
$k=l+1$. From \eqref{s120} we know that for generic
$z'=(z_2,\cdots,z_m)$, we have
\begin{align}\label{s140}
\begin{split}
&\int_0^\delta p(r_1,\cdots,r_m)\frac{r_{k+1}\cdots r_m}{r_{1}
(\log\frac{1}{r_1})^2\cdots r_{k}(\log\frac{1}{r_k})^2}\ dr_1\\
\geq &
\int_0^\delta p(\delta,r_{2},\cdots,r_m)
\frac{r_{k+1}\cdots r_m}{r_{1}
(\log\frac{1}{r_1})^2\cdots r_{k}(\log\frac{1}{r_k})^2}\ dr_1\\
& +\int_0^\delta (\log\log\frac{1}{\delta}-\log\log\frac{1}{r_1})
\frac{r_{k+1}\cdots r_m}{r_{1}
(\log\frac{1}{r_1})^2\cdots r_{k}(\log\frac{1}{r_k})^2}\ dr_1.
\end{split}
\end{align}
The second term in the above formula is integrable with respect
to $r_2,\cdots,r_m$ on $[0,\delta]^{m-1}$ by direct computation.
To estimate the first term, we can choose a generic $\delta$ and
working on the lower dimensional piece $\{ |z_1|=\delta \} \cap
U$. By induction assumption, the first term is also integrable
with respect to $r_2,\cdots,r_m$ on $[0,\delta]^{m-1}$. This
finishes the proof.

\qed

Now we go back to the proof of the main theorem. From the above
lemma we know that $\log f$ is integrable with respect to the
Poincar\'e metric on each chart that touches the divisor $Y$. On
those charts which do not touch $Y$, $\log f $ is bounded. Using
a partition of unity, we can easily see that $\log f \in
L^1( M, \omega_P)$ where $\omega_P$ is the (global)
asymptotic Poincar\'e metric on $ M$. From \eqref{s40}
we have
\begin{align}\label{s150}
\begin{split}
&\int_{ M}\log f
(\bb\partial\bar\partial\rho_{\epsilon})\wedge
(Ric(\omega_{WP})+2m\omega_{WP})^{j-1}\wedge
\widetilde{\omega}^{s-j}\wedge \omega_{0}^{m-s}\\
=&\int_{\rm supp(1-\rho_{\epsilon})}\log f
(\bb\partial\bar\partial\rho_{\epsilon})\wedge
(Ric(\omega_{WP})+2m\omega_{WP})^{j-1}\wedge
\widetilde{\omega}^{s-j}\wedge \omega_{0}^{m-s}\\
=&\sum_{\alpha}\int_{\rm supp(1-\rho_{\epsilon})
\cap U_{\alpha}}\psi_{\alpha}\log f
(\bb\partial\bar\partial\rho_{\epsilon})\wedge
(Ric(\omega_{WP})+2m\omega_{WP})^{j-1}\wedge
\widetilde{\omega}^{s-j}\wedge \omega_{0}^{m-s}
\end{split}
\end{align}
where the sum over $\alpha$ is a finite sum. Since on each
$U_\alpha$, $\log f$ is bounded above by
$c\sum_1^{l_\alpha}\log\log\frac{1}{r_j}$ and both
$Ric(\omega_{WP})+2m\omega_{WP}$ and $\widetilde{\omega}$ are
bounded above and below by a constant multiple of the Poincar\'e
metric, using the same technique as in \eqref{m80}, since the
measure of ${\rm supp}(1-\rho_{\epsilon})$ goes to
$0$ as $\epsilon \to 0$, we conclude that
\[
\int_{ M}\log f
(\bb\partial\bar\partial\rho_{\epsilon})\wedge
(Ric(\omega_{WP})+2m\omega_{WP})^{j-1}\wedge
\widetilde{\omega}^{s-j}\wedge \omega_{0}^{m-s}=0
\]
which implies
\begin{eqnarray}\label{s160}
\int_{ M}
\rho_{\epsilon}(Ric(\omega_{WP})+2m\omega_{WP})^s\wedge
\omega_{0}^{m-s}=\int_{ M}
\rho_{\epsilon}\widetilde{\omega}^s
\wedge \omega_{0}^{m-s}.
\end{eqnarray}
Again, since
$
-c\omega_{P} \leq Ric(\omega_{WP})+2m\omega_{WP} \leq c\omega_{P}
$ and $-c\omega_{P} \leq \widetilde{\omega} \leq c\omega_{P}$,
by the dominate convergence theorem, let $\epsilon \to 0$ in
\eqref{s160} we have
\begin{eqnarray}\label{s170}
\int_{ M}
(Ric(\omega_{WP})+2m\omega_{WP})^s\wedge
\omega_{0}^{m-s}=\int_{ M}
\widetilde{\omega}^s
\wedge \omega_{0}^{m-s}.
\end{eqnarray}
Since
\[
\int_{ M}
\widetilde{\omega}^s
\wedge \omega_{0}^{m-s}=\int_{\bar{ M}}
\widetilde{\omega}^s \wedge \omega_{0}^{m-s}
\]
and $\widetilde{\omega}^s \wedge \omega_{0}^{m-s}$ is a
characteristic class on $\bar{ M}$, we know that
$\int_{ M}\widetilde{\omega}^s
\wedge \omega_{0}^{m-s} \in \mathbb{Z}$ which implies that
$\int_{ M}(Ric(\omega_{WP}))^s\wedge
\omega_{WP}^{m-s} \in \mathbb{Z}$, if the nilpotent operators
are unipotent. In general, by Lemma~\ref{lem41}, the integral
is a rational number.

\qed

We end this paper by the obvious possible generalization
of this paper:

\begin{qu}
Let $c_k(\omega_{WP})$ be the $k$-th elementary polynomial
of the curvature tensor of the Weil-Petersson metric. Let
$X$ be any Weil-Petersson subvariety of dimension
$q$ of a
Weil-Petersson  variety $M$. Then
\[
\int_Xc_k(\omega_{WP})\wedge\omega_{WP}^l
\]
is a rational number, where $k+l=q$. 
\end{qu}

It would be interesting to see if this is true in the
category of Weil-Petersson geometry.

\bibliographystyle{abbrv}   
\bibliography{new041222,unp,local}

\end{document}